\documentclass[11pt,reqno]{amsart}
\usepackage{amsmath, amsfonts, amsthm, amssymb, graphicx}

    \newcommand{\C}{\mathcal{C}}
    
    \newcommand{\ff}{\mathbf{f}}
    \newcommand{\ee}{\mathbf{e}}
    
    \newcommand{\A}{\mathbf{A}}

    \newcommand{\q}{\mathbf{q}}
    \newcommand{\x}{\mathbf{x}}
    \newcommand{\y}{\mathbf{y}}
    
    \newcommand{\rr}{\mathbf{r}}
    \newcommand{\el}{\mathbf{l}}
 \newcommand{\ww}{\mathbf{w}}

\newcommand{\btau}{{\boldsymbol \tau}}

    \newcommand{\del}{{\partial}}

\newcommand{\eps}{\varepsilon}
\newcommand{\cR}{\mathbb{R}}
\newcommand{\fer}[1]{(\ref{#1})}

\newcommand{\R}{\mathbb{R}}

\newcommand{\PP}{\mathcal{P}}

\newcommand {\blp} {\Bigl(}
\newcommand {\brp} {\Bigr)}

\newcommand {\dv} { {\rm div} }

\newcommand{\e} {\varepsilon}

\newcommand{\Sg}{\beta}
\newcommand{\p}{\partial}
\newcommand{\Chi}{\chi}

\newcommand{\n}{{\bf n}}
\renewcommand{\O}{\Omega}
\newcommand{\G}{\Gamma}

\newcommand{\ve}{\varepsilon}
\newcommand{\xxi}{\boldsymbol \xi}
\newcommand{\lam}{\boldsymbol \lambda}

\newcommand{\PtUpL}{{P_1}}
\newcommand{\PtLwL}{{P_2}}
\newcommand{\PtLwR}{{P_3}}
\newcommand{\PtUpR}{{P_4}}

\newtheorem{theorem}{Theorem}[section]

\numberwithin{equation}{section}

\begin{document}
\title[Degenerate Partial Differential Equations]
{On Degenerate Partial Differential Equations}
\author{Gui-Qiang G. Chen}
\address{Gui-Qiang G. Chen, Mathematical Institute, University of Oxford,
         Oxford, OX1 3LB, UK; and Department of Mathematics, Northwestern University,
         Evanston, IL 60208, USA}
\email{\tt chengq@maths.ox.ac.uk;gqchen@math.northwestern.edu}

\date{\today}
\thanks{}

\begin{abstract}
Some of recent developments, including recent results, ideas,
techniques, and approaches, in the study of degenerate partial
differential equations are surveyed and analyzed. Several examples
of nonlinear degenerate, even mixed, partial differential equations,
are presented, which arise naturally in some longstanding,
fundamental problems in fluid mechanics and differential geometry.
The solution to these fundamental problems greatly requires a deep
understanding of nonlinear degenerate partial differential
equations. Our emphasis is on exploring and/or developing unified
mathematical approaches, as well as new ideas and techniques. The
potential approaches we have identified and/or developed through
these examples include kinetic approaches, free boundary approaches, weak
convergence approaches, and related nonlinear ideas and techniques.
We remark that most of the important problems for nonlinear
degenerate partial differential equations are truly challenging and
still widely open, which require further new ideas, techniques, and
approaches, and deserve our special attention and further efforts.
\end{abstract}
\maketitle

\section{Introduction}
We survey and analyze some of recent developments,
including recent results, ideas, techniques, and approaches, in the
study of degenerate partial differential equations. We start with
several important examples of degenerate/mixed linear degenerate
equations and some of their interrelations. Then we present several
examples of nonlinear degenerate, even mixed, partial differential
equations, arising naturally in some longstanding, fundamental
problems in fluid mechanics and differential geometry. These
examples indicate that some of important nonlinear degenerate
problems are ready to be tractable. Our emphasis is on exploring
and/or developing unified mathematical approaches, as well as new
ideas and techniques. The potential approaches we have
identified and/or developed
through these examples include kinetic approaches, free boundary
approaches, weak convergence approaches, and related nonlinear ideas
and techniques.

In fact, nonlinear degenerate, even mixed, partial differential
equations arise also naturally in fundamental problems in many other
areas such as elasticity, relativity, optimization, dynamical
systems, complex analysis, and string theory. The solution to these
fundamental problems in the areas greatly requires a deep
understanding of nonlinear degenerate partial differential
equations. On the other hand,  most of the important problems for
nonlinear degenerate partial differential equations are truly
challenging and still widely open, which require further new ideas,
techniques, and approaches, and deserve our special attention and
further efforts.

During the last half century, three different types of nonlinear
partial differential equations (elliptic, hyperbolic, parabolic)
have been systematically studied separately, and great progress has
been made through the efforts of several generations of
mathematicians. With these achievements, it is the time to
revisit and attack nonlinear degenerate, even mixed, partial
differential equations with emphasis on exploring and developing
unified mathematical approaches, as well as new ideas and techniques.

The organization of this paper is as follows. In Section 2, we
present several important examples of linear degenerate, even mixed,
equations and exhibit some of their interrelations. In Section 3, we
first reveal a natural connection between degenerate hyperbolic
systems of conservation laws and the Euler-Poisson-Darboux equation
through the entropy and Young measure, and then we discuss how this
connection can be applied to solving hyperbolic systems of
conservation laws with parabolic or hyperbolic degeneracy. In
Section 4, we present a kinetic approach to handle a class of
nonlinear degenerate parabolic-hyperbolic equations, the anisotropic
degenerate
diffusion-advection equation. In Section 5, we present two
approaches through several examples to handle nonlinear mixed
problems, especially nonlinear degenerate elliptic problems:
free-boundary techniques and weak convergence methods. In Section 6,
we discuss how the singular limits to nonlinear degenerate
hyperbolic systems of conservation laws via weak convergence methods
can be achieved
through an important limit problem: the vanishing viscosity limit
problem for the Navier-Stokes equations to the isentropic Euler
equations.

\section{Linear Degenerate Equations}
In this section, we present several important examples of linear
degenerate, even mixed, equations and exhibit some of their
interrelations.

Three of the most basic types of partial differential equations are
elliptic, parabolic, and hyperbolic. Consider the partial
differential equations of second-order:
\begin{equation}\label{2.1}
\sum_{i,j=1}^d a_{ij}(\x) \partial_{x_ix_j}u+\sum_{j=1}^d
b_j(\x)\partial_{x_j}u+c(\x)u=f(\x), \quad \x=(x_1,x_2,\cdots,
x_d)\in\Omega\subset\R^d,
\end{equation}
where $a_{ij}(\x), b_j(\x), c(\x)$, and $f(\x)$ are bounded for
$\x\in\Omega$. Equation \eqref{2.1} is called a uniformly elliptic
equation in $\Omega$ provided that there exists $\lambda_0>0$ such that, for any
$\xxi=(\xi_1, \cdots, \xi_d)\in\R^d$,
\begin{equation}\label{2.1-1}
\sum_{i,j=1}^d a_{ij}(\x)\xi_i\xi_j\ge \lambda_0 |\xxi|^2
\qquad\,\, \mbox{for}\,\, \x\in\Omega,
\end{equation}
that is, the $d\times d$ matrix $(a_{ij}(\x))$ is positive definite.

Two of the basic types of time-dependent partial differential
equations of second-order are hyperbolic equations:
\begin{equation}\label{2.2}
\del_{tt} u- L_\x u=f(t,\x),    \qquad t>0,
\end{equation}
and parabolic equations:
\begin{equation}\label{2.3}
\del_{t}u- L_\x u=f(t,\x), \qquad t>0,
\end{equation}
where $L_\x$ is a second-order elliptic operator:
\begin{equation}\label{2.2-1}
L_\x u=\sum_{i,j=1}^d a_{ij}(t,\x)\partial_{x_ix_j}u+\sum_{j=1}^d
b_j(t,\x)\partial_{x_j}u+c(t,\x)u,
\end{equation}
for which $a_{ij}(t,\x), b_j(t, \x), c(t, \x)$, and $f(t,\x)$ are
locally bounded for $(t,\x)\in \R_+^{d+1}:=[0,\infty)\times \R^d$.
Equations \eqref{2.2} and \eqref{2.3} are called uniformly hyperbolic
or parabolic equations, respectively, in a domain in $\R_+^{d+1}$
under consideration, provided that there exists $\lambda_0>0$ such
that, for any $\xxi\in\R^d$,
\begin{equation}\label{2.1-2}
\sum_{i,j=1}^d a_{ij}(t,\x)\xi_i\xi_j\ge \lambda_0 |\xxi|^2
\end{equation}
on the domain,  that is, the $d\times d$ matrix $(a_{ij}(t,\x))$ is
positive definite. In particular, when the coefficient functions
$a_{ij}, b_j, c$, and the nonhomogeneous term $f$ are time-invariant, and
the solution $u$ is also time-invariant, then equations \eqref{2.2}
and \eqref{2.3} are coincide with equation \eqref{2.1}.

Their representatives are Laplace's equation:
\begin{equation}\label{2.7a}
\Delta_\x u=0,
\end{equation}
the wave equation:
\begin{equation}\label{2.9a}
\partial_{tt}u-\Delta_\x u=0,
\end{equation}
and the heat equation:
\begin{equation}\label{2.8a}
\del_t u-\Delta_\x u=0,
\end{equation}
respectively, where $\Delta_\x=\sum_{j=1}^d \partial_{x_jx_j}$ is the Laplace operator.

Similarly, a system of partial differential equations of first-order
in one-dimension:
\begin{equation}\label{2.4}
\del_tU +\A(t,x)\del_x U=0, \qquad U\in\R^n
\end{equation}
is called a strictly hyperbolic system, provided that the $n\times n$
matrix $\A(t,\x)$ has $n$ real, distinct eigenvalues:
$$
\lambda_1(t,x)<\lambda_2(t,x)<\cdots < \lambda_n(t,x)
$$
in the domain under consideration. The simplest example is
\begin{equation}\label{2.10a}
\A(t,x)=\left( \begin{array}{ccc}
0 & -1 \\
-1 & 0\\
\end{array} \right).
\end{equation}
Then system \eqref{2.4}--\eqref{2.10a} with $U=(u,v)^\top$ is
equivalent to the two wave equations:
$$
\partial_{tt}u-\partial_{xx}u=0, \qquad \partial_{tt}v-\partial_{xx} v=0.
$$

However, many important partial differential equations
are degenerate or mixed. That is,
for the linear case, the matrix $(a_{ij}(t,\x))$ or $(a_{ij}(\x))$
is not positive definite or even indefinite, and the eigenvalues of
the matrix $\A(t,x)$ are not distinct.

\subsection{Degenerate Equations and Mixed Equations}
Two prototypes of linear degenerate equations are the {\it
Euler-Poisson-Darboux equation}:
\begin{equation}\label{2.9}
(x-y)\partial_{xy}u+\alpha
\partial_x u+ \beta\partial_y u=0,
\end{equation}
or
\begin{equation}\label{2.10}
x(\partial_{xx}u-
\partial_{yy}u)+ \alpha \partial_x u+ \beta \partial_y u=0,
\end{equation}
and the {\it Beltrami equation}:
\begin{equation}\label{2.8}
 x(\partial_{xx}u+
\partial_{yy}u)+ \alpha\partial_x u+ \beta\partial_y u=0.
\end{equation}

Three prototypes of mixed hyperbolic-elliptic equations are the {\it
Lavrentyev-Betsadze equation}:
\begin{equation}\label{2.11}
\partial_{xx}u+ \mbox{sign}(x)\partial_{yy}u=0
\end{equation}
which exhibits the jump from the hyperbolic phase $x<0$ to the
elliptic phase $x>0$, the {\it Tricomi equation}:
\begin{equation}\label{2.12}
\partial_{xx}u + x\partial_{yy}u=0
\end{equation}
which exhibits hyperbolic degeneracy at $x=0$
(i.e., two eigenvalues coincide at $x=0$,
but the corresponding characteristic curves are not tangential
  to the line $x=0$), and the {\it
Keldysh equation}:
\begin{equation}\label{2.13}
x\partial_{xx}u +\partial_{yy}u=0
\end{equation}
which exhibits parabolic degeneracy at $x=0$ (i.e.,
two eigenvalues coincide at $x=0$,
but the corresponding characteristic curves are
tangential to the line $x=0$).

The mixed parabolic-hyperbolic equations include the linear
degenerate diffusion-advection equation:
\begin{equation}\label{2.14}
\partial_t u+{\mathbf b}\cdot\nabla
u= \nabla\cdot(\A\,\nabla u), \qquad  \x\in \R^d,\,\,
\,t\in\R_+:=(0,\infty),
\end{equation}
where $u: \R_+\times\R^d\to\R$ is unknown, $\nabla=(\partial_{x_1},
\dots,
\partial_{x_d})$ is the gradient operator with respect
to $\x=(x_1, \dots, x_d)\in \R^d$, and ${\mathbf b}:
\R_+\times\R^d\to\R^d$ and $\A: \R_+\times\R^d\to \R^{d\times d}$
are given functions such that $\A=(a_{ij}(t,\x))$ is a $d\times d$
nonnegative, symmetric matrix. When the diffusion matrix function
$\A$ degenerates, the advection term ${\mathbf b}\cdot\nabla u$
dominates; otherwise, the parabolic diffusion term
$\nabla\cdot(\A\,\nabla u)$ dominates.

\subsection{Interrelations between the Linear Equations}

The above linear equations are not actually independent, but are
closely interrelated.

\subsubsection{The Wave Equation \eqref{2.9a}
and the Euler-Poisson-Darboux Equation \eqref{2.10}}

Seek spherically symmetric solutions of the wave equation
\eqref{2.9a}:
$$
v(t,r)=u(t,\x), \qquad r=|\x|.
$$
Then $v(t,r)$ is governed by the Euler-Poisson-Darboux equation
\eqref{2.10} with $\alpha=d-1$ and $\beta=0$:
$$
\partial_{tt}v-\partial_{rr}v -\frac{d-1}{r}\partial_r v=0.
$$
As is well-known, the Euler-Poisson-Darboux equation plays an
important role in the spherical mean method for the wave equation to
obtain the explicit representation of the solution in $\R^d, d\ge 2$
(cf. Evans \cite{Evans-book}).

\subsubsection{The Tricomi Equation \eqref{2.12}, the Beltrami Equation \eqref{2.8},
and the Euler-Poisson-Darboux Equation \eqref{2.10}} Under the
coordinate transformations:
\begin{eqnarray*}
&(x,y)\longrightarrow
(\tau, y)=(\frac{2}{3}(-x)^{\frac{3}{2}}, y),\qquad &\mbox{when}\,\, x<0,\\
&(x,y)\longrightarrow (\tau,
   y)=(\frac{2}{3} x^{\frac{3}{2}}, y),\quad\,\, \qquad &\mbox{when}\,\, x>0,
\end{eqnarray*}
the Tricomi equation \eqref{2.12} is transformed into the Beltrami
equation \eqref{2.8} with $\alpha=\frac{1}{3}$ and $\beta=0$ when $x=\tau>0$:
$$
\partial_{\tau\tau}u +\partial_{yy}u+\frac{1}{3\tau}\partial_\tau
u=0,
$$
and the Euler-Poisson-Darboux equation \eqref{2.10} with
$\alpha=\frac{1}{3}$ and $\beta=0$  when $x=\tau<0$:
$$
\partial_{\tau\tau}u-\partial_{yy}u+\frac{1}{3\tau}\partial_\tau u=0.
$$
These show that a solution to one of them implies the solution of
the other correspondingly, which are equivalent correspondingly.

\medskip
Linear degenerate partial differential equations have been
relatively better understood since 1950. The study of nonlinear
partial differential equations has been focused mainly on the
equations of single type during the last half century. The three
different types of nonlinear partial differential equations have
been systematically studied separately, and one of the main focuses
has been on the tools, techniques, and approaches to understand different
properties and features of solutions of the equations with these
three different types. Great progress has been made through the
efforts of several generations of mathematicians.

\medskip
As we will see through several examples below, nonlinear degenerate,
even mixed, partial differential equations naturally arise in some
fundamental problems in fluid mechanics and differential geometry.
The examples include nonlinear degenerate hyperbolic systems of
conservation laws, the nonlinear degenerate diffusion-advection
equation and the Euler equations for compressible flow in fluid
mechanics, and the Gauss-Codazzi system for isometric embedding in
differential geometry. Such degenerate, or mixed, equations
naturally arise also in fundamental problems in many other areas
including elasticity, relativity, optimization, dynamical systems,
complex analysis, and string theory. The solution to these
fundamental problems in the areas greatly requires a deep
understanding of nonlinear degenerate, even mixed, partial differential
equations.

\section{Nonlinear Degenerate Hyperbolic Systems of Conservation Laws}

Nonlinear hyperbolic systems of conservation laws in one-dimension take the
following form:
\begin{equation}
\del_t U + \del_x F(U) = 0, \,\qquad U \in \R^n, \,\, (t,x)\in
\R_+^2. \label{2.1a}
\end{equation}
For any $C^1$ solutions, \eqref{2.1a} is equivalent to
$$
\del_t U +\nabla F(U) \del_x U = 0, \,\qquad U \in \R^n, \,\,
(t,x)\in \R_+^2.
$$
Such a system is hyperbolic if the $n\times n$ matrix $\nabla F(U)$
has $n$ real eigenvalues $\lambda_j(U)$ and linearly independent
eigenvectors $\rr_j(U), 1\le j\le n$. Denote
\begin{equation}
\mathcal{D}:=\{U: \lambda_i(U)=\lambda_j(U), \quad i\ne j, \, 1\le
i, j\le n\} \label{2.2a}
\end{equation}
as the degenerate set. If the set $\mathcal{D}$ is empty, then this
system is strictly hyperbolic and, otherwise, degenerate hyperbolic.
Such a set allows a degree of interaction, or nonlinear resonance,
among different characteristic modes, which is missing in the
strictly hyperbolic case but causes additional analytic
difficulties. A point $U_* \in \mathcal{D}$ is hyperbolic degenerate
if $\nabla F(U_*)$ is diagonalizable and, otherwise, is parabolic
degenerate.

Degenerate hyperbolic systems of conservation laws have arisen from
many important fields such as continuum mechanics including the
vacuum problem, multiphase flows in porous media, MHD, and
elasticity. On the other hand, degenerate hyperbolicity of systems
is generic in some sense. For example, for three-dimensional
hyperbolic systems of conservation laws, Lax \cite{Lax3} indicated
that systems with $2( mod\, 4)$ equations must be degenerate
hyperbolic. The result is also true when the systems have $\pm 2,
\pm 3, \pm 4 (mod\, 8)$ equations (see \cite{FRS}). Then the plane
wave solutions of such systems are governed by the corresponding
one-dimensional hyperbolic systems with $\mathcal{D}\ne \emptyset$.

Since $F$ is a nonlinear function, solutions of the Cauchy problem
for \eqref{2.1a} (even starting from smooth initial data) generally
develop singularities in a finite time, and then the solutions
become discontinuous functions. This situation reflects in part the
physical phenomenon of breaking of waves and development of shock
waves. For this reason, attention is focused on solutions in the space
of discontinuous functions, where one can not directly use the
classical analytic techniques that predominate in the theory of
partial differential equations of other types.

To overcome this difficulty, a natural idea is to construct
approximate solutions $U^\epsilon(t,x)$ to \eqref{2.1a} by using
shock capturing methods and then to study the convergence and
consistency of the approximate solutions to \eqref{2.1a}. The key
issue is whether the approximate solutions converge in an
appropriate topology and the limit function is  consistent with
\eqref{2.1a}. Solving this issue involves two aspects: one is to
construct good approximate solutions, and the other is to make
suitable compactness analysis in an appropriate topology.

\subsection{Connection with the Euler-Poisson-Darboux Equation:
Entropy and Young Measure}
\medskip

The connection between degenerate hyperbolic systems of conservation
laws and the Euler-Poisson-Darboux equation is through the entropy.

 A pair of functions $(\eta, q): \R^n  \to \R^2 $ is called
an entropy-entropy flux pair (entropy pair, in short) if they
satisfy the following linear hyperbolic system:
\begin{equation}\label{3.1a}
\nabla q(U) = \nabla\eta(U) \nabla F(U).
\end{equation}
Then the function $\eta(U)$ is called an entropy. Clearly,
any $C^1$ solution satisfies
\begin{equation}\label{3.2a} \del_t \eta(U)
+ \nabla\eta(U)\nabla F(U)\del_x U=0,
\end{equation}
or
\begin{equation}\label{3.3a}
\del_t \eta(U) + \del_x q(U)  = 0.
\end{equation}
For a $BV$ solution that is not $C^1$, the second term in
\eqref{3.2a} has no meaning in the classical sense because of the
multiplication of a Radon measure with a discontinuous function. If
$\eta(U)$ is an entropy, then the left side of \eqref{3.2a} becomes
the left side of \eqref{3.3a} that makes sense even for $L^\infty$
solutions in the sense of distributions. An $L^\infty$ function
$U(t,x)$ is called an entropy solution (cf. Lax \cite{Lax2}) if
\begin{equation}
\del_t \eta(U) + \del_x q(U) \le 0 \label{3.4a}
\end{equation}
in the sense of distributions for any convex entropy pair, that is,
the Hessian $\nabla^2\eta(U)\ge 0$.

Assume that system \eqref{2.1a} is endowed with globally defined
Riemann invariants $\ww(U)=(w_1, \dots, w_n)(U), 1\le i \le n$, which satisfy
\begin{equation}
\nabla w_i(U)\cdot \nabla F(U)=\lambda_i(U)\nabla w_i(U).
\label{3.5a}
\end{equation}
 The necessary and sufficient condition
for the existence of the Riemann invariants $w_i(U), 1\le i\le n$, for strictly
hyperbolic systems is the well-known Frobenius condition:
$$
\el_i\{\rr_j, \rr_k\}=0, \qquad \text{for any}\,\,\, j, k \ne i,
$$
where $\el_i$ denotes the left eigenvector corresponding to
$\lambda_i$ and $\{\cdot, \cdot\}$ is the Poisson bracket of vector
fields in the $U$-space (cf. \cite{AMa,Law}).
For $n=2$, the Frobenius condition always
holds. For any smooth solution $U(t,x)$, the corresponding Riemann
invariants $w_i, 1\le i\le n,$ satisfy the transport equations
$$
\del_t w_i + \lambda_i(\ww) \del_x w_i = 0,
$$
which indicate that $w_i$ is invariant along the $i$-th
characteristic field.

Taking the inner product of \eqref{3.5a} with the right eigenvectors
$\rr_i$ of $\nabla F$ produces the characteristic form
$$
(\lambda_i\nabla\eta-\nabla q)\cdot \rr_i =0,
$$
that is,
\begin{equation}\label{3.6a}
\lambda_i \partial_{w_i}\eta =\partial_{w_i} q.
\end{equation}

For $n=2$, the linear hyperbolic system \eqref{3.5a} is equivalent
to \eqref{3.6a}, and system \eqref{3.6a} can be reduced to the
following linear second-order hyperbolic equation:
\begin{equation}
\partial_{w_1w_2}\eta+\frac{\partial_{w_1}\lambda_{2}}{\lambda_2-\lambda_1}\partial_{w_2}\eta
-\frac{\partial_{w_2}\lambda_{1}}{\lambda_2-\lambda_1}\partial_{w_1}\eta=0.
\label{3.7a}
\end{equation}
For the case $\mathcal{D}=\emptyset$, equation \eqref{3.7a} is, in
general, regular and hyperbolic for which either the Goursat problem
or the Cauchy problem is well posed in the coordinates of Riemann
invariants. The entropy space is infinite-dimensional and is
represented by two families of functions of one variable. However,
for the case $\mathcal{D}\ne \emptyset$, the situation is much more
complicated because of the singularity of the functions
$\frac{\partial_{w_1}\lambda_{2}}{\lambda_2-\lambda_1}$ and
$\frac{\partial_{w_2}\lambda_{1}}{\lambda_2-\lambda_1}$ on the set
$\ww(\mathcal{D})\subset\R^2$. The typical form
of such equations is
\begin{equation}
\partial_{w_1w_2}\eta +\frac{\alpha(\frac{w_1}{w_2},
w_1-w_2)}{w_1-w_2}\partial_{w_1}\eta +\frac{\beta(\frac{w_1}{w_2},
w_1-w_2)}{w_1-w_2}\partial_{w_2}\eta=0. \label{3.8a}
\end{equation}
This equation may have extra singularities of $\alpha$ and $\beta$
both at the origin and on the line $w_1=w_2$. The questions are
whether there exist nontrivial regular solutions of the singular
equation \eqref{3.8a} and, if so, how large the set of smooth
regular solutions is.

The connection between the compactness problem of approximate
solutions and the entropy determined by the Euler-Poisson-Darboux
equation \eqref{3.7a} is the Young measure via the compensated
compactness ideas first developed by Tartar \cite{Ta1,Ta2} and Murat
\cite{Mu1,Mu2} and a related observation presented by Ball
\cite{Ball}.

The Young measure is a useful tool for studying the limiting
behavior of measurable function sequences. For an arbitrary sequence
of measurable maps
$$
U^\epsilon: \R^2_+ \to \R^n, \quad \|U^\epsilon\|_{L^\infty}\le C
<\infty,
$$
that converges in the weak-star topology of $L^\infty$ to a function
$U$,
$$
w^*-\lim\limits_{\epsilon \to 0} U^\epsilon = U,
$$
there exist a subsequence (still labeled) $U^\epsilon$ and a family
of Young measures
$$
\nu_{t,x}\in \text{Prob}.(\R^n), \quad \text{supp}~\nu_{t,x} \subset
\{\lam: |\lam|\le M\}, \qquad (t,x)\in \R^2_+,
$$
such that, for any continuous function $g$,
$$
w^*-\lim\limits_{\epsilon\to 0}g(U^\epsilon(t,x))
=\int_{\R^n}g(\lam)d\nu_{t,x}(\lam):=\langle\nu_{t,x}, g\rangle
$$
for almost all points $(t,x)\in \R^2_+$. In particular,
$U^\epsilon(t,x)$ converges strongly to $U(t,x)$ if and only if the
Young measure $\nu_{t,x}$ is equal  to a Dirac mass concentrated at
$U(t,x)$ for a.e. $(t,x)$.

In many cases, one can estimate (cf. \cite{DiPerna1,Chen2}) that the
approximate solutions $U^\epsilon(t,x)$ generated by the shock
capturing methods for \eqref{3.1a} satisfy

\begin{itemize}
\item  $\|U^\epsilon\|_{L^\infty}\le C < \infty$.
\item For $C^2$ entropy pairs $(\eta_i, q_i), i=1,2,$
determined by \eqref{3.6a} and \eqref{3.7a},
\begin{equation}
\del_t
\eta_i(U^\epsilon) + \del_x q_i(U^\epsilon) \quad \text{compact in}
\, \, H^{-1}_{\text{loc}}. \label{3.9a}
\end{equation}
\end{itemize}
Then,
for any $C^2$ entropy pairs $(\eta_i, q_i), i=1,2$, determined by
\eqref{3.6a} and \eqref{3.7a}, the Young measure $\nu_{t,x}$ is
forced to satisfy the Tartar-Murat commutator relation:
\begin{equation}\label{3.11a}
\langle\nu_{t,x},
    \eta_1 q_2-q_1\eta_2\rangle
=\langle \nu_{t,x}, \eta_1 \rangle \langle\nu_{t,x}, q_2\rangle
 -\langle\nu_{t,x}, q_1 \rangle\langle\nu_{t,x}, \eta_2 \rangle
\qquad\,\, a.e. \,\, (t,x)\in\R^2_+.
\end{equation}

\medskip
\noindent
This relation is derived from \eqref{3.9a}, the Young
measure representation theorem for the measurable function sequence,
and a basic continuity theorem for the $2\times 2$ determinant in
the weak topology (cf. \cite{Ta1,Chen2}). This indicates that
proving the compactness of the approximate solutions
$U^\epsilon(t,x)$ is equivalent to solving the functional equation
\eqref{3.11a} for the Young measure $\nu_{t,x}$ for all possible
$C^2$ entropy pairs $(\eta_i, q_i), i=1,2$, determined by
\eqref{3.6a} and \eqref{3.7a}. If one clarifies the structure of the
Young measure $\nu_{t,x}(\lam)$, one can understand the limiting
behavior of the approximate solutions. For example, if one can prove
that almost all Young measures $\nu_{t,x}, (t,x)\in \R^2_+$, are
Dirac masses, then one can conclude the strong convergence of the
approximate solutions $U^\epsilon(t,x)$ almost everywhere.
On the other hand, the structure of the Young measures is determined by the
$C^2$ solutions of the Euler-Poisson-Darboux equation \eqref{3.7a}
via \eqref{3.11a}.
One of the principal difficulties for the reduction is the general
lack of enough classes of entropy functions that can be verified to
satisfy certain weak compactness conditions. This is due to possible
singularities of entropy functions near the regions with degenerate
hyperbolicity. The larger the set of $C^2$ solutions to the entropy
equation \eqref{3.7a} of the Euler-Poisson-Darboux type is, the
easier solving the functional equation \eqref{3.11a} for the Young
measure is.

\subsection{Hyperbolic Conservation Laws with Parabolic Degeneracy}
One of the prototypes for such systems is the system of isentropic
Euler equations with the form
\begin{equation}
\begin{cases}
\del_t\rho + \del_x m =0, \\
\del_t m + \del_x(\frac{m^2}{\rho}+p(\rho))= 0,
\end{cases}
\label{4.1a}
\end{equation} where $\rho, m$, and $p$ are the
density, mass, and pressure, respectively. For $\rho>0$, $v=m/\rho$
represents the velocity of the fluid. The physical region for
\eqref{4.1a} is $\{(\rho, m)\, :\, |m|\leq C \, \rho\}$ for some
$C>0$, in which the term $\frac{m^2}{\rho}$ in the flux function is
only at most Lipschitz continuous near the vacuum. For $\rho>0$,
$v=\frac{m}{\rho}$ represents the velocity of the fluid.

The pressure $p$ is a function of the density through the internal
energy $e(\rho)$:
\begin{equation}\label{3.15a}
p(\rho) =\rho^2\, e'(\rho) \qquad \mbox{for}\,\, \rho\ge 0.
\end{equation}
In particular, for a polytropic perfect gas,
\begin{equation}\label{Eq:Pressure}
p(\rho)=\kappa\rho^\gamma, \qquad e(\rho)=\int_0^\rho
\frac{p(r)}{r^2}dr=\frac{\kappa}{\gamma-1}\rho^{\gamma-1},
\end{equation}
where $\gamma>1$ is the adiabatic exponent and, by the scaling, the
constant $\kappa$ in the pressure-density relation may be chosen as
$\kappa=\frac{(\gamma-1)^2}{4\gamma}$ without loss of generality.

For \eqref{4.1a}, strict hyperbolicity and genuine nonlinearity away
from the vacuum require that
\begin{equation}\label{IEE-1}
p'(\rho)>0, \quad 2 \, p'(\rho)+\rho \, p''(\rho)>0 \qquad
\text{for} \,\,\,\, \rho>0.
\end{equation}
Near the vacuum, for some $\gamma>1$,
\begin{equation}\label{IEE-2}
\frac{p(\rho)}{\rho^\gamma}\,\to\, \kappa_1>0 \qquad\mbox{when
$\rho\to 0$}.
\end{equation}

One of the fundamental features of this system is that strict
hyperbolicity fails when $\rho\to 0$. That is,
$\mathcal{D}=\{\rho=0\}$, which is the vacuum state, and the
degeneracy is parabolic.

The eigenvalues of \eqref{4.1a} are
$$
\lambda_1=\frac{m}{\rho}-c(\rho), \qquad
\lambda_2=\frac{m}{\rho}+c(\rho),
$$
and the corresponding Riemann invariants of \eqref{4.1a} are
$$
w_1=\frac{m}{\rho}-k(\rho), \qquad w_2=\frac{m}{\rho}+k(\rho),
$$
where
$$
c(\rho):=\sqrt{p'(\rho)}
$$
is the sound speed, and
$$
k(\rho):= \int_0^\rho \frac{c(r)}{r} \, dr.
$$
Then the entropy equation \eqref{3.7a} in the Riemann coordinates is
the following Euler-Poisson-Darboux equation:
\begin{equation}
\partial_{w_1w_2}\eta
+\frac{\alpha(w_1-w_2)}{w_1-w_2}(\partial_{w_1}\eta-\partial_{w_2}\eta)=0,
\label{4.2}
\end{equation}
where
$$
\alpha(w_1-w_2) =-\frac{k(\rho)k''(\rho)}{k'(\rho)^2}
$$
with $\rho=-k^{-1}(\frac{w_1-w_2}{2})$.

For the $\gamma$-law gas,
$$
\alpha(w_1-w_2) =\lambda:=\frac{3-\gamma}{2(\gamma-1)}
$$
is a constant, the simplest case.

In terms of the variables $(\rho, v)$, $\eta$ is determined by the
following second-order linear wave equation:
\begin{equation}
\partial_{\rho\rho}\eta - k'(\rho)^2 \,\partial_{vv} \eta = 0. \label{4.2a}
\end{equation}
An entropy $\eta$ is called a weak entropy if
$\eta(0,v) = 0$.
For example, the mechanical energy
pair:
$$
\eta_\ast=\frac{1}{2}\frac{m^2}{\rho}+\rho e(\rho), \quad q_*
=\frac{m^3}{2\rho^2}+ m\big(e(\rho)+\frac{p(\rho)}{\rho}\big),
$$
is a convex weak entropy pair.

By definition, the weak {\it entropy kernel\/} is the solution
$\chi(\rho,v,s)$ of the singular Cauchy problem
\begin{equation}\label{4.2b}
\begin{cases}
\partial_{\rho\rho}\chi- k'(\rho)^2 \,\partial_{vv} \chi = 0,\\
\chi(0,v,s) = 0, \\
\partial_\rho \chi(0,v,s) = \delta_{v=s}
\end{cases}
\end{equation}
in the sense of distributions, where $s$ plays the role of a
parameter and $\delta_{v=s}$ denotes the Dirac measure at $v=s$.
Then the family of weak entropy functions is described by
\begin{equation}\label{4.2c}
\eta(\rho,v) = \int_\R \chi(\rho, v,s) \,\psi(s) \, ds,
\end{equation}
where $\psi(v)$ is an arbitrary function. By construction,
$\eta(0,v) = 0$, $\eta_\rho(0,v) = \psi(v)$.
One can prove that, for $0\le\rho\le
C$, $|\frac{m}{\rho}|\le C$,
$$
|\nabla\eta(\rho, m)|\le C_\eta, \qquad
 |\nabla^2\eta(\rho, m)|\le C_\eta\nabla^2\eta_\ast(\rho, m),
$$
for any weak entropy $\eta$, with $C_\eta$ independent of $(\rho,
m)$.

Since this system is a prototype in mathematical fluid dynamics, the
mathematical study of this system has an extensive history dating back to the
work of Riemann \cite{Ri}, where a special Cauchy problem, so-called
Riemann problem, was solved. Zhang-Guo \cite{ZG} established an
existence theorem of global solutions to this system for a class
of initial value functions by using the characteristic method.
Nishida \cite{Ni} obtained the first large data existence theorem
with locally finite total-variation for the case $\gamma=1$ using
Glimm's scheme \cite{Glimm1}. Large-data theorems have also been
obtained for general $\gamma >1$ in the case where the initial value
functions with locally finite total-variation are restricted to
prevent the development of cavities (e.g. \cite{DZWH,Liu1,NS}). The
difficult point in bounding the total-variation norm at low
densities is that the coupling between the characteristic fields
increases as the density decreases. This difficulty is a reflection
of the fact that the strict hyperbolicity fails at the vacuum:
$\mathcal{D}\ne\emptyset$.

The first global existence for \eqref{4.1a} with large initial data
in $L^\infty$ was established in DiPerna \cite{DiPerna1} for the
case $\gamma=\frac{N+2}{N}$, $N\ge 5$ odd, by the vanishing
viscosity method. The existence problem for usual gases with general
values $\gamma \in (1,\frac{5}{3}]$ was solved in Chen \cite{Chen1} and
Ding-Chen-Luo \cite{DCL} (also see \cite{Chen2}). The case
$\gamma\ge 3$ was treated by Lions-Perthame-Tadmor \cite{LPT}.
Lions-Perthame-Souganidis \cite{LPS} dealt with the interval $(\frac{5}{3},
3)$ and simplified the proof for the whole interval. A compactness
framework has been established even for the general pressure law in
Chen-LeFloch \cite{ChenLeFloch} by using only weak entropy pairs.

More precisely, assume the pressure function $p=p(\rho)\in
C^4(0,\infty)$ satisfies condition \eqref{IEE-1} (i.e., strict
hyperbolicity and genuine nonlinearity) away from the vacuum and,
near the vacuum, $p(\rho)$ is only asymptotic to the $\gamma$-law
pressure (as real gases): there exist a sequence of exponents
\begin{equation}\label{IEE-pa}
1<\gamma:=\gamma_1<\gamma_2< \ldots < \gamma_N \leq (3\gamma -1)/2 <
\gamma_{N+1}
\end{equation}
and a sufficiently smooth function $P=P(\rho)$ such that
\begin{eqnarray}
&& p(\rho) = \sum_{j=1}^N \kappa_j \, \rho^{\gamma_j} +
\rho^{\gamma_{N+1}} \, P(\rho), \label{IEE-pb}
\\
&& P(\rho) \, \text{ and } \, \rho^3 \, P'''(\rho) \, \text{ are
bounded as } \rho \to 0,
\label{IEE-pc}
\end{eqnarray}
for some coefficients $\kappa_j \in \R$ with $\kappa_1>0$. The
solutions under consideration will remain in a bounded subset of
$\{\rho\ge 0\}$ so that the behavior of $p(\rho)$ for large $\rho$
is irrelevant. This means that the pressure law $p(\rho)$ has the
same singularity as $\sum_{j=1}^N \kappa_j \, \rho^{\gamma_j}$ near
the vacuum. Observe that $p(0)=p'(0)=0$, but, for $k>\gamma_1$, the
higher derivative $p^{(k)}(\rho)$ is unbounded near the vacuum with
different orders of singularity.

\begin{theorem}[Chen-LeFloch \cite{ChenLeFloch}]
Consider system  \eqref{4.1a} with general pressure law satisfying
\eqref{IEE-1} and \eqref{IEE-pa}--\eqref{IEE-pc}. Assume that a
sequence of functions $(\rho^\ve, m^\ve)$ satisfies that

{\rm (i)} There exists $C>0$ independent of $\ve$ such that
\begin{equation}\label{3.8}
0\leq\rho^\ve(t,x)\leq C, \qquad |m^\ve(t,x)|\leq C \, \rho^\ve(t,x)
\qquad \hbox{ for a.e. } (t,x);
\end{equation}

{\rm (ii)} For any weak entropy pair $(\eta, q)$ of \eqref{4.1a},
\eqref{IEE-1}, and \eqref{IEE-pa}--\eqref{IEE-pc},
\begin{equation}\label{3.9}
 \del_t\eta(\rho^\ve, m^\ve)
+\del_xq(\rho^\ve,m^\ve) \quad \hbox{ is compact in }
H^{-1}_{loc}(\R_+^2).
\end{equation}
 Then the sequence $(\rho^\ve,
m^\ve)$ is compact in $L^1_{loc}(\R_+^2)$. Moreover, there exists a
global entropy solution $(\rho(t,x), m(t,x))$ of the Cauchy problem
\eqref{4.1a}, \eqref{IEE-1}, and \eqref{IEE-pa}--\eqref{IEE-pc},
satisfying
$$
0\le\rho(t,x)\le C, \qquad |m(t,x)|\le C\rho(t,x)
$$
for some $C$ depending only on $C_0$ and $\gamma$, and
$$
\del_t\eta(\rho,m) + \del_x q(\rho,m) \le 0
$$
in the sense of distributions for any {\it convex} weak entropy pair
$(\eta, q)$. Furthermore, the bounded solution operator
$(\rho,m)(t,\cdot)=S_t(\rho_0,m_0)(\cdot)$ is compact in $L^1$ for
$t>0$.
\end{theorem}

As discussed in Section 3.1, if the Young measure satisfying
\eqref{3.11a} reduces to a Dirac mass for a.e.~$(t,x)$, then the
sequence $(\rho^\varepsilon, m^\varepsilon)$ is compact in the
strong topology and converges subsequentially toward an entropy
solution. For the Euler equations, to obtain that the Young measure
$\nu_{(t,x)}$ is a Dirac mass in the $(\rho,m)$-plane, it suffices
to prove that the measure in the $(\rho,v)$-plane, still denoted by
$\nu_{(t,x)}$, is either a single point or a subset of the vacuum
line
$$
\{(\rho,v)\, :
\,\rho=0,\,|v|\leq\sup_{\eps>0}\big\|\frac{m^\eps}{\rho^\eps}\big\|_{L^\infty}\}.
$$
The main difficulty is that only {\it weak\/} entropy pairs can be
used, because of the vacuum problem.

In the proof of \cite{Chen1,DCL,DiPerna1} (also cf. ~\cite{Chen2}),
the heart of the matter is to construct special functions $\psi$ in
\eqref{4.2c} in order to exploit the form of the set of constraints
\eqref{3.11a}. These test-functions are suitable approximations of
high-order derivatives of the Dirac measure. It is used that
\eqref{3.11a} represents an {\it imbalance of regularity:\/} the
operator on the left is more regular than the one on the right due
to cancellation. DiPerna \cite{DiPerna1} considered the case that
$\lambda:=\frac{3-\gamma}{2(\gamma-1)}\ge 3 $ is an integer so that
the weak entropies are polynomial functions of the Riemann
invariants. The new idea of applying the technique of fractional
derivatives was first introduced in Chen \cite{Chen1} and
Ding-Chen-Luo \cite{DCL} to deal with real values of $\lambda$.

A new technique for equation \eqref{3.11a} was introduced by
Lions-Perthame-Tadmor \cite{LPT} for $\gamma \in [3, \infty)$ and by
Lions-Perthame-Souganidis \cite{LPS} for $\gamma \in (1,3)$.
Motivated by a kinetic formulation, they made the
observation that the use of the test-functions $\psi$ could in fact
be bypassed, and \eqref{3.11a} be directly expressed with the
entropy kernel $\chi_*$. Namely, \eqref{3.11a} holds for all $s_1$
and $s_2$ by replacing $\eta_j:=\chi_*(s_j)$ and $q_j
:=\sigma_*(s_j)$ for $j=1,2$. Here $\sigma_*$ is the entropy-flux
kernel defined as
$$
\sigma_*(\rho,v,s)
         = \big(v + \theta \,(s-v)\big)\, \chi_*(\rho,v,s).
$$
In \cite{LPS}, the commutator relations were differentiated in $s$,
using the technique of fractional derivatives as introduced in
\cite{Chen1,DCL} by performing the operator $\del_s^{\lambda+1}$, so
that singularities arise by differentiation of $\chi_*$. This
approach relies again on the lack of balance in regularity of the
two sides of \eqref{3.11a} and on the observation that
$\langle\nu_{(t,x)},\chi_*(s)\rangle$ is smoother than the kernel
$\chi_*(s)$ itself, due to the averaging by the Young measure.

However, many of the previous arguments do not carry over to the general
pressure law. The main issue is to construct all of the weak entropy
pairs of \eqref{4.1a}.
The proof in Chen-LeFloch \cite{ChenLeFloch} has been based on the
existence and regularity of the {\it entropy kernel\/} that
generates the family of weak entropies via solving a new
Euler-Poisson-Darboux equation, which is {\it highly singular\/}
when the density of the fluid vanishes. New properties of {\it
cancellation of singularities\/} in combinations of the entropy
kernel and the associated entropy-flux kernel have been found. In
particular,
a new {\it multiple-term\/} expansion has been introduced based on
the Bessel functions with suitable exponents, and the optimal
assumption required on the singular behavior on the pressure law at
the vacuum has been identified in order to valid the multiple-term
expansion and to establish the existence theory. The results cover,
as a special example, the density-pressure law $p(\rho) = \kappa_1
\, \rho^{\gamma_1} + \kappa_2 \, \rho^{\gamma_2}$ where $\gamma_1,
\gamma_2 \in (1,3)$ and $\kappa_1, \kappa_2 >0$ are arbitrary
constants.
The proof of the reduction theorem for Young measure has also
further simplified the proof known for the polytropic perfect gas.

Then this compactness framework has been successfully applied to
proving the convergence of the Lax-Friedrichs scheme, the Godunov
scheme, and the artificial viscosity method for the isentropic
Euler equations with the general pressure law.

\subsection{Hyperbolic Conservation Laws with Hyperbolic Degeneracy}

One of the prototypes of hyperbolic conservation laws with
hyperbolic degeneracy is the gradient quadratic flux system, which
is umbilic degeneracy, the most singular case:
\begin{equation}
\del_t U + \del_x (\nabla_U C(U)) = 0, \quad U=(u,v)^\top \in \R^2,
\label{4.3a}
\end{equation}
 where
\begin{equation} C(U)= \frac{1}{3}au^3+ bu^2v + uv^2, \label{4.4a}
\end{equation}
and $a$ and $b$ are two real parameters.

Such systems are generic in the following sense. For any smooth
nonlinear flux function, take its Taylor expansion about the
isolated umibilic point. The first three terms including the
quadratic terms determine the local behavior of the hyperbolic
singularity near the umbilic point. The hyperbolic degeneracy
enables us to eliminate the linear term by a coordinate
transformation to obtain the  system with a homogeneous quadratic
polynomial flux. Such a polynomial flux contains some inessential
scaling parameters. There is a nonsingular linear coordinate
transformation to transform the above system into \eqref{4.3a}--
\eqref{4.4a}, first studied by Marchesin, Isaacson, Plohr, and
Temple, and in a more satisfactory form by Schaeffer-Shearer
\cite{SS1}. {}From the viewpoint of group theory, such a reduction
from six to two parameters is natural: For the six dimensional space
of quadratic mappings acted by the four dimensional group $GL(2,
\R)$, one expects the generic orbit to have codimension two.

The Riemann solutions for such systems were discussed by Isaacson,
Marchesin, Paes-Leme, Plohr, Schaeffer, Shearer, Temple, and others
(cf. \cite{IMPT,IT,SS2,SSMP}). Two new types of shock waves, the
overcompressive shock  and the undercompressive shock, were
discovered, which are quite different from the gas dynamical shock.
The overcompressive shock can be easily understood using the Lax
entropy condition \cite{Lax2}. It is known that there is a traveling
wave solution connecting the left to right state of the
undercompressive shock for the artificial viscosity system.
Stability of such traveling waves for the overcompressive shock and
the undercompressive shock has been studied (cf. \cite{LK,LX}).

The next issues are whether the compactness of the corresponding
approximate solutions is affected by the viscosity matrix as the
viscosity parameter goes to zero, to understand the sensitivity of
the undercompressive shock with respect to the viscosity matrix, and
whether the corresponding global existence of entropy solutions can
be obtained as a corollary from this effort. The global existence of
entropy solutions to the Cauchy problem for a special case of such
quadratic flux systems was solved by Kan \cite{Ka} via the viscosity
method. A different proof  was given independently to the same
problem by Lu \cite{Lu}.

In Chen-Kan \cite{ChenKan}, an $L^\infty$ compactness framework has
been established for sequences of approximate solutions to general
hyperbolic systems with umbilic degeneracy specially including
\eqref{4.3a}--\eqref{4.4a}. Under this framework, any approximate
solution sequence, which is apriori uniformly bounded in $L^\infty$
and produces the correct entropy dissipation, leads to the
compactness of the corresponding Riemann invariant sequence. This
means that the viscosity matrix does not affect the compactness of
the corresponding uniformly bounded Riemann invariant sequence.
Again, one of the principal difficulties associated with such
systems is the general lack of enough classes of entropy functions
that can be verified to satisfy certain weak compactness conditions.
This is due to possible singularities of entropy functions near the
regions of degenerate hyperbolicity. The analysis leading to the
compactness involves two steps:

In the first step, we have constructed regular entropy functions
governed by a highly singular entropy equation, the
Euler-Poisson-Darboux type equation,
\begin{equation}
\partial_{w_1w_2}\eta+\frac{\alpha(\frac{w_1}{w_2})}{w_2-w_1}\partial_{w_1}\eta
+\frac{\beta({\frac{w_1}{w_2}})}{w_2-w_1}\partial_{w_2}\eta=0. \label{4.5a}
\end{equation}
There are two main difficulties. The first is that the coefficients
of the entropy equation are multiple-valued functions near the
umbilic points in the Riemann invariant coordinates. This difficulty
has been overcome by a detailed analysis of the singularities of the
Riemann function of the entropy equation. This analysis involves a
study of a corresponding Euler-Poisson-Darboux equation and requires
very complicated estimates and calculations. An appropriate choice
of Goursat data leads to the cancellation of singularities and the
achievement of regular entropies in the Riemann invariant
coordinates. The second difficulty is that the nonlinear
correspondence between the physical state coordinates and the Riemann
invariant coordinates is, in general, irregular. A regular entropy
function in the Riemann invariant coordinates is usually no longer
regular in the physical coordinates. We have overcome this by a
detailed analysis of the correspondence between these two
coordinates.

In the second step, we have studied the structure of the Young
measure associated with the approximate sequence and have proved
that the support of the Young measures lies in finite isolated
points or separate lines in the Riemann invariant coordinates. This
has been achieved by a delicate use of Serre's technique
\cite{Serre} and regular entropy functions, constructed in the first
step, in the Tartar-Murat commutator equation for Young measures
associated with the approximate solution sequence.

This compactness framework has been successfully applied to proving
the convergence of the Lax-Friedrichs scheme, the Godunov scheme,
and the viscosity method for the quadratic flux systems. Some
corresponding existence theorems of global entropy solutions for
such systems have been established. The compactness has been
achieved by reducing the support of the corresponding Young measures
to a Dirac mass in the physical space.

\section{Nonlinear Degenerate Parabolic-Hyperbolic Equations}

One of the most important examples of nonlinear degenerate
parabolic-hyperbolic equations is the nonlinear degenerate
diffusion-advection equation of second-order with the form:
\begin{equation}
\p_t u+\nabla\cdot \ff(u)=\nabla\cdot (\A(u)\nabla u), \qquad \x\in
\R^d, \quad t\in\R_+, \label{PSCL1}
\end{equation}
where $\ff: \cR\to \cR^d$ satisfies $\ff^\prime (\cdot) \in
L^\infty_{\rm loc}(\cR; \cR^d)$, and the $d\times d$ matrix $\A(u)$
is symmetric, nonnegative, and locally bounded, so that
$\A(u)=(\sigma_{ik}(u))(\sigma_{ik}(u))^\top$ with
$\sigma_{ik}(u)\in L^\infty_{loc}(\cR)$. Equation \eqref{PSCL1} and
its variants model anisotropic degenerate diffusion-advection
motions of ideal fluids and arise in a wide variety of important
applications, including two-phase flows in porous media and
sedimentation-consolidation processes (cf. \cite{BCBT,CJ,EFM,Vaz}).
Because of its importance in applications, there is a large
literature for the design and analysis of various numerical methods
to calculate the solutions of equation \eqref{PSCL1} and its
variants (e.g. \cite{CJ,CoD,DDE,EFM,KR1}) for which a well-posedness
theory is in great demand.

\medskip
One of the prototypes is the porous medium equation:
\begin{equation}\label{PorousMedium}
\p_t u=\Delta_\x (u^m), \qquad m> 1,
\end{equation}
which describes the fluid flow through porous media (cf.
\cite{Vaz}). Equation \eqref{PorousMedium} is degenerate on the
level set $\{u=0\}$; away from this set, i.e., on $\{u>0\}$, the
equation is strictly parabolic. Although equation
\eqref{PorousMedium} is of parabolic nature, the solutions exhibit
certain hyperbolic feature, which results from the degeneracy. One
striking family of solutions is Barenblatt's solutions found in
\cite{Barenblatt}:
$$
u(t,\x; a, \tau)=\frac{1}{(t+\tau)^k}
\Big[a^2-\frac{k(m-1)}{2md}\frac{|\x|^2}{(t+\tau)^{\frac{2k}{d}}}\Big]_+^{\frac{1}{m-1}},
$$
where $k=\frac{d}{(m-1)d+2}$, and $a\ne 0$ and $\tau>0$ are any
constants.
The dynamic boundary of the support of $u(t,\x;a,\tau)$ is
$$
|\x|=\sqrt{\frac{2m d}{k(m-1)}}a (t+\tau)^{\frac{k}{d}}, \qquad t\ge 0.
$$
{\it This shows that the support of Barenblatt's
solutions propagates with a finite speed!}

The simplest example for the isotropic case (i.e., $\A(u)$ is
diagonal) with both phases is
\begin{equation}\label{Isotropic}
\p_t u+ \p_x(\frac{u^2}{2})=\p_{xx} [u]_+,
\end{equation}
where $[u]_+=\max\{u, 0\}$ (cf. \cite{ChenDiB}). Equation
\eqref{Isotropic} is hyperbolic when $u<0$ and parabolic when $u>0$,
and the level set $\{u=0\}$ is a free boundary that is an interface
separating the hyperbolic phase from the parabolic phase. For any
constant states $u^\pm$ such that $u^+<0<u^-$ and $y_*\in\R$, there
exists a unique nonincreasing profile $\phi=\phi(y)$ such that
\begin{eqnarray*}
&&\lim_{y\to \pm \infty}\phi(y)=u^\pm, \qquad \lim_{y\to
  \pm\infty}\phi'(y)=0,\\
&&\phi(y)\equiv u^+ \,\qquad\text{for
  all}\,\, y>y_*,\\
&&\phi\in C^2(-\infty, y_*),\,\,\,  \phi(y_*-)=0,\,\,\,
   \lim_{y\to y_*-}\phi'(y)=-\infty, \\
&& \frac{d[\phi(y)]_+}{dy}\Big|_{y=y_*-}
   =\frac{1}{2}u^+u^-
\end{eqnarray*}
so that $u(t,x)=\phi(x-st)$ is a discontinuous solution connecting
$u^+$ from $u^-$, where $s=\frac{u^++u^-}{2}$. Although $u(t,x)$ is
discontinuous, $[u(t,x)]_+$ is a continuous function even across the
interface $\{u=0\}$.
See Fig. 1.

\begin{figure}[h]
\centering
\includegraphics[scale=0.7]{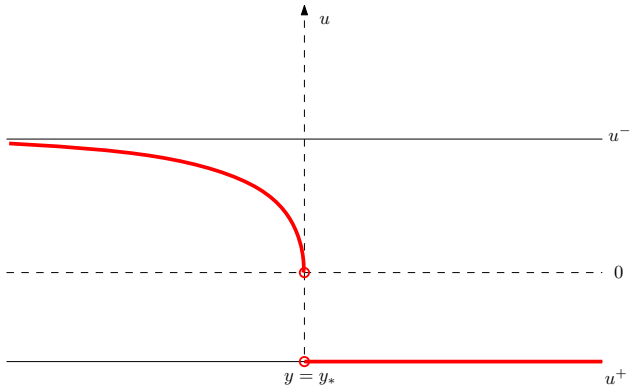}
\caption[]{A discontinuous profile connecting the hyperbolic phase
from the parabolic phase} \label{fig:SRDP}
\end{figure}

\smallskip
The well-posedness issue for the Cauchy problem is relatively well
understood if one removes the diffusion term
$\nabla\cdot(\A(u)\nabla u)$, thereby obtaining a scalar hyperbolic
conservation law; see Lax \cite{Lax1}, Oleinik \cite{Oleinik},
Volpert \cite{Volpert}, Kruzhkov \cite{Kr}, and
Lions-Perthame-Tadmor \cite{LPT}, and Perthame \cite{Pe1}
(also cf. \cite{Da,Sm}). It is equally well understood for the
diffusion-dominated case, especially when the set $\{u :
rank(\A(u))<d\}$ consists of only isolated points with certain order
of degeneracy; see Brezis-Crandall \cite{BC}, Caffarelli-Friedman
\cite{CaF},
 Daskalopoulos-Hamilton \cite{DH},
 DiBenedetto \cite{Di}, Gilding \cite{Gilding}, V\'{a}zquez \cite{Vaz},
and the references cited therein. For the isotropic diffusion,
$a_{ij}(u)=0, i\ne j$, some stability results for entropy solutions
were obtained for $BV$ solutions by Volpert-Hudjaev \cite{VH} in
1969. Only in 1999, Carrillo \cite{Carrillo} extended this result to
$L^\infty$ solutions (also see Eymard-Gallou\"{e}t-Herbin-Michel
\cite{EGHM}, Karlsen-Risebro \cite{KR1} for further extensions), and
Chen-DiBenedetto \cite{ChenDiB} handled the case of unbounded
entropy solutions which may grow when $|\x|$ is large.

A unified approach, the kinetic approach, to deal with both
parabolic and hyperbolic phases for the general anisotropic case
with solutions in $L^1$ has been developed in Chen-Perthame
\cite{ChenPerthame1}. This approach is motivated by the macroscopic
closure procedure of the Boltzmann equation in kinetic theory, the
hydrodynamic limit of large particle systems in statistical
mechanics, and early works on kinetic schemes to calculate shock
waves and theoretical kinetic formulation for the pure hyperbolic
case; see
\cite{Br,Cercignani,ChapmanCowling,kipnis,LPT,OVY,Pe1,Pe2,Sone} and
the references cited therein. In particular, a notion of kinetic
solutions and a corresponding kinetic formulation have been
extended. More precisely, consider the kinetic function,
quasi-Maxwellian, $\Chi$ on $\cR^2$:
\begin{equation}
     \Chi(v;u)=\left \{
      \begin{array}{ll}
      +1 &\ \ \,\, {\rm for}\ \ 0 < v < u, \\
      -1 &\ \ \,\, {\rm for}\ \ u < v < 0,\\
      \;0 &\ \ \,\, {\rm otherwise}.
      \end{array}
      \right.
\label{PDefChi}
\end{equation}
If $u\in L^\infty\big(\R_+; L^1(\cR^d)\big)$, then $\Chi(v;u)\in
L^\infty\big(\R_+; L^1(\cR^{d+1})\big)$.

\medskip
\noindent {\bf Definition}. A function $u(t,\x)\in
L^\infty\big(\R_+; L^1(\cR^d)\big)$ is called a kinetic solution if
$u(t,\x)$ satisfies the following:

\smallskip
\noindent{\rm (i)} The kinetic equation:
\begin{equation}\label{KS1}
\p_t\Chi(v;u) +\ff'(v)\cdot \nabla \Chi(v;u) -\nabla\cdot
(\A(v)\nabla\Chi(v;u)) =\p_v (m+ n)(t,\x,v)
\end{equation}
holds in the sense of distributions with the initial data $
\Chi(v;u)|_{t=0} = \Chi(v;u_0), $ for some nonnegative measures
$m(t,\x;v)$ and $n(t,\x;v)$, where $n(t,\x,v)$ is defined by
\begin{equation}
\langle\, n(t,\x,\cdot), \psi(\cdot)\rangle:=\sum\limits_{k=1}^{d}
   \blp \sum\limits_{i=1}^{d}\p_{x_i}\Sg^\psi_{ik}(u(t,\x))\brp^2
   \in L^1(\R_+^{d+1}),
\label{KS3}
\end{equation}
for any $\psi\in C_0^\infty(\cR)$ with $\psi \geq 0$ and $
\beta_{ik}^{\psi}(u):=\int^u\sigma_{ik}(v)\sqrt{\psi(v)}\,dv; $

\noindent {\rm (ii)} There exists $\mu\in L^\infty(\R)$ with $0\le
\mu(v)\to 0$ as $|v|\to\infty$ such that
\begin{equation}\label{NoSupport}
\int_0^\infty \int_{\cR^d} (m+n)(t,\x; v)\,dt\, d\x \leq\mu(v);
\end{equation}

\noindent{\rm (iii)} For any two nonnegative functions $\psi_1,\,
\psi_2\in C_0^\infty(\cR)$,
\begin{equation}
\sqrt{\psi_1(u(t,\x))}
 \sum\limits_{i=1}^{d}\p_{x_i}\Sg^{\psi_2}_{ik}(u(t,\x))
= \sum\limits_{i=1}^{d}\p_{x_i}\Sg^{\psi_1 \psi_2}_{ik}(u(t,\x))
\qquad a.e.\,\,\, (t,\x). \label{KS2}
\end{equation}

\medskip
Then we have

\begin{itemize}
\item {\bf Well-posedness in $L^1$}: Under this notion, the space
$L^1$ is both a well-posed space for kinetic solutions and a
well-defined space for the kinetic equation in (i). That is, the
advantage of this notion is that the kinetic equation is well
defined even when the macroscopic fluxes and diffusion matrices are
not locally integrable so that $L^1$ is a natural space on which the
kinetic solutions are posed. This notion also covers the so-called
renormalized solutions used in the context of scalar hyperbolic
conservation laws by B{\'e}nilan-Carrillo-Wittbold \cite{BCW}. Based
on this notion, a new approach has been developed in
\cite{ChenPerthame1} to establish the contraction property of
kinetic solutions in $L^1$. This leads to a well-posedness
theory---existence, stability, and uniqueness---for the Cauchy
problem for kinetic solutions in $L^1$.

\medskip
\item {\bf Consistency of the kinetic equation with
the original macroscopic equation}: When the kinetic solution $u$ is
in $L^\infty$, for any $\eta\in C^2$ with $\eta''(u)\ge 0$,
multiplying $\eta'(v)$ both sides of the kinetic equation in (i) and
then integrating in $v\in\R$ yield
$$
\partial_t\eta(u)+\nabla_\x\cdot(\q(u)-\A(u)\nabla_\x\eta(u))
=-\int_{\R}\eta''(v)(m+n)(t,\x;v)dv\le 0.
$$
In particular, taking $\eta(u)=\pm u$ yields that $u$ is a weak
solution to the macroscopic equation. The  uniqueness result implies
that any kinetic solution in $L^\infty$ must be an entropy solution.
On the other hand, any entropy solution is actually a kinetic
solution. Therefore, the two notions are equivalent for solutions in
$L^\infty$, although the notion of kinetic solutions is more
general.

\medskip
\item For the isotropic case, condition (iii) automatically
holds, which is actually a chain rule. In fact, the extension from
the isotropic to anisotropic case is not a purely technical issue,
and the fundamental and natural chain-rule type property (iii) does
not appear in the isotropic case and turns out to be the
corner-stone for the uniqueness in the anisotropic case. Moreover,
condition (ii) implies that $m+n$ has no support  at $u=\infty$.
\end{itemize}

\bigskip
Based on this notion, the corresponding kinetic formulation, and the
uniqueness proof in the pure hyperbolic case in \cite{Pe1}, we have
developed a new effective approach to establish the contraction
property of kinetic solutions in $L^1$. This leads to a
well-posedness theory for the Cauchy problem of \eqref{4.1a} with
initial data:
\begin{equation}\label{PSCL2}
u|_{t=0}=u_0 \in L^\infty(\cR^d)
\end{equation}
for kinetic solutions only in $L^1$.

\begin{theorem}[Chen-Perthame \cite{ChenPerthame1}]\label{Thm:Main}
\noindent {\rm (i)} For any kinetic solution $u\in L^\infty (\R_+;
L^1(\cR^d))$ with initial data $u_0(\x)$, we have
$$
\|u(t,\cdot)-u_0 \|_{L^1(\cR^d)} \to 0 \qquad\text{as}\,\,\, t \to
0;
$$

\noindent {\rm (ii)} If $u, v\in L^\infty (\R_+;L^1(\cR^d))$ are
kinetic solutions to {\rm \fer{PSCL1}} and {\rm \fer{PSCL2}} with
initial data $u_0(\x)$ and $v_0(\x)$, respectively, then
\begin{equation}\label{contractionp}
\|u(t,\cdot)-v(t, \cdot) \|_{L^1(\cR^d)} \leq \|u_0 -v_0
\|_{L^1(\cR^d)};
\end{equation}

\noindent{\rm (iii)} For initial data $u_0\in L^1(\cR^d)$, there
exists a unique kinetic solution $u\in C\big(\R_+; L^1(\cR^d)\big)$
for the Cauchy problem \eqref{4.1a} and \eqref{PSCL2}. If $u_0\in
L^\infty\cap L^1(\cR^d)$, then the kinetic solution is the unique
entropy solution and $\|u(t,\cdot)\|_{L^\infty(\R^d)} \leq \| u_0
\|_{L^\infty(\cR^d)}$.

Furthermore, assume that the flux function $\ff(u)\in C^1$ and the
diffusion matrix $\A(u)\in C$ satisfy the nonlinearity-diffusivity
condition: The set
\begin{equation}\label{nonlinearity}
\{v\,:\, \tau +\ff'(v)\cdot \y=0, \,\,
                 \y \A(u)\,\y^\top=0\}\subset \R
\end{equation}
has zero Lebesgue measure, for any $\tau\in\R, \, \y=(y_1,\dots,
y_d)$ with $|\y|=1$. Let $u\in L^\infty (\R_+^{d+1})$ be the unique
entropy solution to {\rm \fer{4.1a}} and {\rm \fer{PSCL2}} with
periodic initial data $u_0\in L^\infty$ for period
${\mathbb T}_P=\Pi_{i=1}^d[0, P_i]$, i.e., $u_0(\x+P_i\ee_i)=u_0(\x)$ a.e.,
where $\{\ee_i\}_{i=1}^d$ is the standard basis in $\R^d$. Then
\begin{equation}\label{decay}
\big\|u(t,\cdot)-\frac{1}{|{\mathbb T}_P|}\int_{{\mathbb T}_P}
u_0(\x)\,d\x\big\|_{L^1({\mathbb T}_P)} \to 0 \qquad\,\,\,
\text{as}\,\, \, t \to \infty,
\end{equation}
where $|{\mathbb T}_P|$ is the volume of the period ${\mathbb
T}_P$.
\end{theorem}

The {\em nonlinearity-diffusivity condition} implies that there is
no interval of $v$ in which the flux function $\ff(v)$ is affine and
the diffusion matrix $\A(v)$ is degenerate. Unlike the pure
hyperbolic case, equation {\rm \fer{PSCL1}} is no longer
self-similar invariant, and the diffusion term in the equation
significantly affects the behavior of solutions; thus the argument
in Chen-Frid \cite{ChenFrid} based on the self-similar scaling for
the pure hyperbolic case could not be directly applied. The argument
for proving Theorem \ref{Thm:Main} is based on the kinetic approach
developed in \cite{ChenPerthame1}, involves a time-scaling and a
monotonicity-in-time of entropy solutions, and employs the
advantages of the kinetic equation \eqref{KS1}, in order to
recognize the role of nonlinearity-diffusivity of equation
\eqref{PSCL1} (see \cite{ChenPerthame2}).

\smallskip
Based on the very construction of the kinetic approach,
the results
can easily be translated in terms of the old
Kruzkov entropies by integrating in $v$. In the case of uniqueness
for the general case, this was performed in Bendahmane-Karlsen
\cite{BKa}.

\smallskip
Follow-up results based on the Chen-Perthame approach in
\cite{ChenPerthame1} include $L^1$-error estimates and continuous
dependence of solutions in the convection function and the diffusion
matrix in Chen-Karlsen \cite{ChenKarlsen1}; and more general
degenerate {diffusion-advection-reaction equations} in Chen-Karlsen
\cite{ChenKarlsen2}. Other recent developments include the related
notion of dissipative solutions in Perthame-Souganidis \cite{PS},
regularity results of solutions in Tadmor-Tao \cite{TT}, as well as
different types of diffusion terms in Andreianov-Bendahmane-Karlsen
\cite{ABK}.

\section{Nonlinear Equations of Mixed Hyperbolic-Elliptic Type}

Unlike the linear case, very often, nonlinear equations can not be
separated as two degenerate equations, but are truly mixed. Nonlinear
partial differential equations of mixed hyperbolic-elliptic type
arise naturally in many fundamental problems. In this
section, we present two approaches through several examples to
handle nonlinear mixed problems, especially nonlinear degenerate
elliptic problems: Free-boundary techniques and weak convergence
methods.

\subsection{Weak Convergence Methods}

We first present two problems: transonic flow pass obstacles and
isometric embedding, for which weak convergence methods, especially
methods of compensated compactness, play an important role.

\subsubsection{Transonic Flow Pass Obstacles in Gas Dynamics}
By scaling, the Euler equations for compressible, isentropic, irrotational
fluids take the form:
\begin{equation}\label{5.1}
\begin{cases}
\del_x(\rho u)+ \del_y(\rho v) = 0, \\
\del_xv -\del_y u = 0,
\end{cases}
\end{equation}
combined with the Bernoulli relation:
\begin{equation}\label{5.3}
\rho = \Big(1 - \frac{\gamma-1} 2  (u^2+v^2)\Big)^{\frac{1}{\gamma-1}}, \qquad \gamma>1,
\end{equation}
for the pressure-density relation: $p=p(\rho)=\rho^\gamma/\gamma$.
This provides two equations for the two unknowns $(u,v)$.
Furthermore, we note that, if $\rho$ is constant (which is the
incompressible case), the two equations in \eqref{5.1} become the
Cauchy-Riemann equations for which any boundary value problem can be
posed for the elliptic partial differential equations.

By the second equation in \eqref{5.1}, we introduce the velocity
potential $\varphi$:
\begin{equation}\label{5.4}
(u,v) = \nabla \varphi.
\end{equation}
Then
our conservation law of mass becomes a nonlinear partial
differential equation of second-order for $\varphi$:
\begin{equation}\label{5.5}
 \del_x(\rho \,\del_x\varphi) + \del_y(\rho\, \del_y\varphi) = 0,
\end{equation}
which is combined with the Bernoulli relation:
\begin{equation}\label{5.6}
\rho = \Big(1 -{\gamma - \frac{1}{2}} |\nabla
\varphi|^2\Big)^{\frac{1}{\gamma - 1}}.
\end{equation}

Introduce the sound speed $c$ as
\begin{equation}\label{5.7}
 c^2 = p' (\rho ) = 1 - \frac{\gamma - 1}{2} q^2,
 \qquad \mbox{with the fluid speed}\,\, q=|\nabla\varphi|=\sqrt{u^2+v^2},
\end{equation}
so that, at the sonic value when $q = c$, we have $q = q_{*}$ with
\begin{equation}\label{5.8}
q_{*} := \sqrt {\frac{2}{\gamma + 1}} .
\end{equation}
Then equation \eqref{5.5} is elliptic if
$$
q < q_{*}
$$ and hyperbolic
when
$$
q > q_{*}.
$$
There is an upper bound placed on $q$ from the
Bernoulli relation:
\begin{equation}\label{5.9}
 q \leq q_{cav}:= \sqrt{\frac{2}{\gamma-1}},
\end{equation}
where $q_{cav}$ is the cavitation speed for which $\rho = 0$.

On the other hand,
equation \eqref{5.5} corresponds to the Euler-Lagrange equation
for the functional
\begin{equation}\label{5.11}
    \int_\Omega G( | \nabla \varphi |) dx dy,
\end{equation}
where
\begin{equation}\label{5.10}
G(q)=\int^{q^2}\Big(1-{\frac{\gamma-1}{2}}s \Big)^{\frac{1}{\gamma-1}} ds.
\end{equation}
Since $(q_{*}-q)G''(q)>0$,
the direct method of calculus of variations  (e.g. Evans
\cite{Evans-book}) provides the existence of weak solutions if
it is known apriori that the flow is subsonic $(q < q_{*})$ so that
$G$ is convex and the problem is elliptic.  For example, this
includes the fundamental problem of subsonic flow around a profile
as formulated in Bers \cite{Bers-book}.

A profile $\PP$ is prescribed by a smooth curve, except for a
trailing edge with an opening $\varepsilon \pi$ at $z_T$, $0 \leq
\varepsilon \leq 1$.

\begin{figure}[h]
\centering
\includegraphics[height=0.5in,width=4.0in]{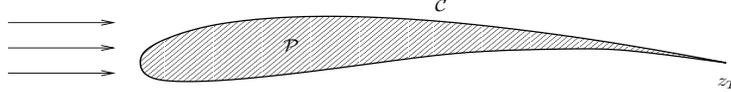}  
\caption[]{Profile $\C$ of the obstacle $\mathcal{P}$: $\e=0$ }
\label{fig:NF1}
\end{figure}

\bigskip
\bigskip
\bigskip

\begin{figure}[h]
\centering
\includegraphics[height=0.5in,width=4.0in]{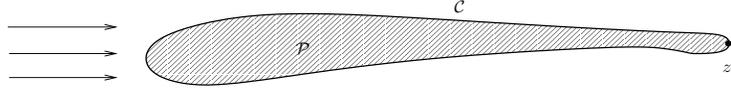}  
\caption[]{Profile $\C$ of the obstacle $\mathcal{P}$: $\e=1$ }
\label{fig:NF2}
\end{figure}

\bigskip
\bigskip
\bigskip

\begin{figure}[h]
\centering
\includegraphics[height=0.5in,width=4.0in]{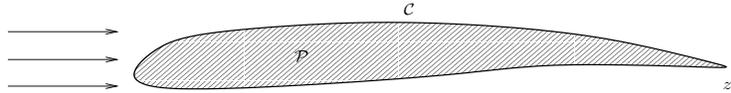}  
\caption[]{Profile $\C$ of the obstacle $\mathcal{P}$: $0<\e<1$ }
\label{fig:NF3}
\end{figure}

If $\varepsilon=0$, the profile has a tangent at the trailing edge
(see Fig. \ref{fig:NF1}).
The tangent to $\PP$ satisfies a uniform H\"older condition with
respect to the arc length.  The velocity $\ww=(u, v)$ must attain a given
subsonic limit at infinity.  We enforce the Kutta-Joukowski
condition:
\begin{eqnarray*}
&q \rightarrow 0 \qquad &\mbox{as $(x,y) \rightarrow z_T$ if
$\varepsilon = 1$},\\
&q = O(1) \qquad  &\mbox{as $(x,y) \rightarrow z_T$ if $0 \leq
\varepsilon < 1$},
\end{eqnarray*}
and define problem $P_1 (\ww_\infty)$ with a
prescribed constant velocity $\ww_\infty$ at infinity.

For a smooth profile, $\varepsilon =1$ (see Fig. \ref{fig:NF2}), define the circulation
\begin{equation}\label{5.12}
    \Gamma = \oint_\PP (u,v) \cdot \btau\, ds,
\end{equation}
where $\btau$ is the unit tangent to $\PP$. In this case, we can
also consider problem $P_2 (\ww_\infty, \Gamma)$ where the circulation
is prescribed, instead of the Kutta-Joukowski condition.

In both problems,  the slip boundary condition on the profile is required:
\begin{equation}\label{5.13}
(u,v)\cdot {\bf n}= 0\qquad\quad \mbox{on} \,\, \PP \qquad \mbox
{(boundary condition)}
\end{equation}
where ${\bf n}$ is the exterior unit normal on $\PP$.

The first existence theorem for $P_1(\ww_\infty)$ was given by
Frankl-Keldysh \cite{FK} for sufficiently small speed at infinity.
For a general gas, the first complete existence theorem for
$P_2(\ww_\infty,\Gamma)$ was given by Shiffman \cite{Shif1}.  This was
followed by a complete existence and uniqueness theorem by Bers
\cite{Bers2} for $P_1(\ww_\infty)$, a stronger uniqueness result of
Finn-Gilbarg \cite{FG3}, and a higher dimensional result of Dong
\cite{Dong-book}.  The basic result is as follows:

\smallskip
{\it For a given constant velocity at infinity, there exists a
number $\hat q$ depending upon the profile $\PP$ and the adiabatic exponent
$\gamma>1$ such that the problem $P_1(\ww_\infty)$ has a unique
solution for $0 < q_\infty:=|\ww_\infty| < \hat q$. The maximum $q_m$
of $q$ takes on all values between $0$ and $q_{*}$, $q_m\rightarrow
0$ as $q_\infty\rightarrow 0$, and $q_m\rightarrow q_{*}$ as
$q_\infty\rightarrow \hat q$. A similar result holds for
$P_2(\ww_\infty, \Gamma)$.}

\smallskip
The main tool for the results is to know apriori that, if $q_\infty
< \hat q$ (i.e., the speed at infinity is not only subsonic) but
sufficiently subsonic, then $q < q_{*}$ in the whole flow domain.
Subsonic flow at infinity itself does not guarantee that the flow
remains subsonic, since the profile produces flow orthogonal to the
original flow direction.
Shiffman's proof did use the direct method of the calculus of
variations, while Bers's relied on both elliptic methods and the
theory of pseudo-analytic functions.  The existence of a critical
point for the variational problem would be a natural goal for the case
when $q_\infty$ is not restricted to be less than $\hat q$,
since it would provide a direct proof of our boundary value problem.
However, no such proof has been given.

More recent investigations based on weak convergence methods start
in the 1980's. DiPerna \cite{DiPerna3} suggested that the
Murat-Tartar method of compensated compactness be amenable to flows
which exhibit both elliptic and hyperbolic regimes, and investigated
an asymptotic approximation to our system called the steady
transonic small disturbance equation (TSD). He proved that, if a list of
assumptions were satisfied (which then guaranteed the applicability
of the method of compensated compactness), then a weak solution
exists to the steady TSD equation.
However, no one has ever been able to show that DiPerna's list
is indeed satisfied.

In \cite{Mor85} (also see
\cite{morawetz-04}), Morawetz layed out a program for proving the
existence of the steady transonic flow problem about a bump profile
in the upper half plane (which is equivalent to a symmetric profile
in the whole plane). In particular, Morawetz showed
that, if the key hypotheses of the method of compensated compactness
could be satisfied, now known as a ``compactness framework'' (see
Chen \cite{GQChen}), then indeed there would exist a weak solution
to the problem of flow over a bump which is exhibited by subsonic
and supersonic regimes, i.e., transonic flow.

The ``compactness framework'' for our system can be stated as follows:
A sequence of functions $\ww^\varepsilon
(x,y)=(u^\varepsilon, v^\varepsilon)(x,y)$ defined on an open set
$\Omega\subset\R^2$ satisfies the following  set of conditions:

\smallskip
\noindent (A.1)\,
 $q^\varepsilon(x,y)=|\ww^\varepsilon (x,y) | \leq q_* \;\; \hbox{a.e. in} \;\; \Omega$
for some positive constant $q_* < q_{cav}$;

\smallskip
\noindent (A.2) \, $\partial_x Q_{1 \pm} (\ww^\varepsilon ) +
\partial_y Q_{2 \pm} (\ww^\varepsilon) \; \text{are confined in a
compact set in } H^{-1}_{\rm loc}(\Omega)$ for entropy pairs
$(Q_{1\pm}, Q_{2\pm})$, and
$(Q_{1 \pm} (\ww^\varepsilon), Q_{2 \pm} (\ww^\varepsilon
))$ are confined to a bounded set uniformly in $L_{\rm loc}^\infty
(\Omega)$,
where $(Q_1, Q_2)$ is an entropy pair, that is,
$\partial_x Q_1 (\ww^\varepsilon ) +
\partial_y Q_2(\ww^\varepsilon ) = 0$ along smooth solutions of our
system).

\smallskip
In case (A.1) and (A.2) hold, then
the Young measure $\nu_{x,y}$ determined by
the uniformly bounded sequence of functions $\ww^\e(x,y)$
is constrained by the following commutator relation:
\begin{equation}\label{5.14}
\begin{split}
   \langle \nu_{x,y},\;  Q_{1+} Q_{2-} - Q_{1-}  Q_{2+}\rangle
    = \langle \nu_{x,y}, \; Q_{1+}\rangle \langle \nu_{x,y}, \; Q_{2-}\rangle
 - \langle \nu_{x,y}, \; Q_{1-} \rangle
\langle \nu_{x,y}, \; Q_{2+}\rangle.
\end{split}
\end{equation}

The main point for the compensated compactness framework is to prove
that $\nu_{x,y}$ is a Dirac measure by using entropy pairs, which
implies the compactness of the sequence
$\ww^\varepsilon(x,y)=(u^\varepsilon, v^\varepsilon)(x,y)$ in
$L^1_{\rm loc}(\Omega)$.  In this context, both DiPerna
\cite{DiPerna3} and Morawetz \cite{Mor85} needed to presume
the existence of an approximating sequence parameterized by
$\varepsilon$ to their problems satisfying (A.1) and (A.2) so that
they could exploit the commutator relation and obtain the strong
convergence in $L^1_{\rm loc}(\Omega)$ to a weak solution of their
problems.

As it turns out, there is one problem where (A.1) and (A.2) hold
trivially, i.e., the sonic limit of subsonic flows.  In that case,
we return to the result by Bers \cite{Bers2} and Shiffman
\cite{Shif1}, which says that, if the speed at infinity $q_\infty$
is less than some $\hat q$, there is a smooth unique solution to
problems $P_1 (\ww_\infty)$ and $P_2 (\ww_\infty, \Gamma)$ and ask what
happens as $q_\infty \nearrow \hat q$.  In this case, the flow
develops sonic points and the governing equations become degenerate
elliptic. Thus, if we set $\varepsilon=\hat q-q_\infty$ and examine
a sequence of exact smooth solutions to our system, we see trivially
that (A.1) is satisfied since $|q^\varepsilon|\leq q_{*}$, and
(A.2) is satisfied since $\partial_x Q_\pm (\ww^\varepsilon)
+\partial_y Q_\pm (\ww^\varepsilon)=0$ along our solution  sequence.
The effort is in finding entropy pairs which can guarantee that the
Young measure $\nu_{x,y}$ reduces to a Dirac mass. Ironically, the
original conservation equations of momentum in fact provide two sets
of entropy pairs, while the irrotationality and mass conservation
equations provide another two sets.  This observation has been
explored in detail in Chen-Dafermos-Slemrod-Wang \cite{CDSW}.

What then about the fully transonic problem of flow past an obstacle
or bump where $q_\infty > \hat q$? In Chen-Slemrod-Wang \cite{CSW1},
we have provided some of the ingredients for satisfying (A.1) and
(A.2). More precisely, we have introduced the usual flow angle
$\theta = \tan^{-1}(\frac{v}{u})$ and written the irrotationality
and mass conservation equation as an artificially viscous problem:
\begin{equation} \label{5.15}
\begin{cases}
    \partial_x v^\e -\partial_y u^\e =\varepsilon\Delta\theta^\e ,  \\
   \partial_x (\rho^\e u^\e) + \partial_y (\rho^\e v^\e)
   =\varepsilon\nabla\cdot(\sigma(\rho^\e)\nabla \rho^\e),
\end{cases}
\end{equation}
where $\sigma (\rho)$ is suitably chosen, and appropriate boundary
conditions are imposed for this regularized ``viscous'' problem. The
crucial point is that a uniformly $L^\infty$ bound in
$q^\varepsilon$ has been obtained when $1\leq\gamma <3$ which uniformly
prevents cavitation.  However, in this formulation,
a uniform bound in the flow angle
$\theta^\varepsilon$ and lower bound of $q^\varepsilon$ away
from zero (stagnation) in any fixed region disjoint
from the profile are still  assumed apriori to guarantee the
$(q,\theta)$--version of (A.1).  On the other hand, (A.2) is easily
obtained from the viscous formulation by using a special entropy
pair of Osher-Hafez-Whitlow \cite{OHW}.  In fact, this entropy pair
is very important:  It guarantees that the inviscid limit of the
above viscous system satisfies a physically meaningful ``entropy''
condition (Theorem 2 of \cite{OHW}). With (A.1)
and (A.2) satisfied,
Morawetz's theory \cite{Mor85} then applies to yield the strong
convergence in $L^1_{\rm loc} (\Omega)$ of our approximating
sequence. It would be interesting to establish
a uniform bound in the flow angle
$\theta^\varepsilon$ and lower bound of $q^\varepsilon$ away
from zero (stagnation) in any fixed region disjoint
from the profile under some natural conditions on the profile.

\subsubsection{Isometric Embeddings in Differential Geometry}
In differential geometry, a longstanding, fundamental problem is to
characterize intrinsic metrics on a two-dimensional surface
$\mathcal{M}^2$ which can be realized as embeddings into $\R^3$. Let
$\O\subset\R^2$ be an open set and $g_{ij}, i, j=1,2,$ be a given
matric on $\mathcal{M}^2$ parameterized on $\Omega$. Then the first
fundamental form for $\mathcal{M}^2$ is
\begin{equation}\label{5.16}
I=g_{11}(dx)^2+2 g_{12}dxdy +g_{22}(dy)^2.
\end{equation}

\medskip
{\bf Isometric Embedding Problem.} {\it Seek an injective map
$\mathbf{r}: \Omega\to \R^3$ such that $d\mathbf{r}\cdot
d\mathbf{r}=I$, i.e.,
\begin{equation}\label{5.17}
\p_x{\bf r}\cdot\p_x{\bf r}=g_{11},\quad
    \p_x{\bf r}\cdot\p_y{\bf r}=g_{12},\quad
   \p_y{\bf r}\cdot\p_y{\bf r}=g_{22},
\end{equation}
so that $(\p_x{\bf r}, \p_y{\bf r})$ in $\R^3$ are linearly
independent.}

\medskip
The equations above are three nonlinear partial differential
equations for the three components of ${\bf r}$. Recall that the
second fundamental form $I\!I$ for $\mathcal{M}^2$ defined on
$\Omega$ is
\begin{equation}\label{5.18}
I\!I=-d\n\cdot d{\bf r}=h_{11}(dx)^2 + 2 h_{12}dxdy +h_{22} (dy)^2,
\end{equation}
and $h_{ij}$ is the orthogonality of the unit normal $\n$ of
the surface ${\bf r}(\Omega)\subset\R^3$ to its tangent plane. The
Christoffel symbols are
$$
\G_{ij}^{(k)}:=\frac{1}{2}g^{kl}\left(\partial_j
g_{il}+\partial_i g_{jl}-\partial_l
 g_{ij}\right),
$$
which depend on the first derivatives of $g_{ij}$, where the
summation convention is used, $(g^{kl})$ denotes the inverse of
$(g_{ij})$, and $(\partial_1, \partial_2):=(\p_x, \p_y)$.

The fundamental theorem of surface theory indicates that there
exists a surface in $\R^3$ whose first and second fundamental
 forms are $I$ and $I\!I$ if the coefficients $(g_{ij})$ and
 $(h_{ij})$ of the two given quadratic forms $I$ and $I\!I$, $I$
 being positive definite, satisfy the Gauss-Codazzi system.
This theorem also holds even for discontinuous coefficients $h_{ij}$
(cf. Mardare \cite{Mardare1}).

Given  $(g_{ij})$, the coefficients  {$(h_{ij})$} of $I\!I$ are
determined by the {\it Gauss-Codazzi system}, which consists of the
Codazzi equations:
\begin{equation}\label{5.19}
\left\{\begin{aligned}
&\del_x{M}-\del_y{L}=\G^{(2)}_{22}L-2\G^{(2)}_{12}M+\G^{(2)}_{11}N, \\
&\del_x{N}-\del_y{M}=-\G^{(1)}_{22}L +2\G^{(1)}_{12}M-\G^{(1)}_{11}N,
\end{aligned}\right.
\end{equation}
and the Gauss equation, i.e., {\it Monge-Amp\`{e}re} constraint:
\begin{equation}\label{5.20}
LN-M^2=K,
\end{equation}
where
\begin{equation}\label{5.21}
L=\frac{h_{11}}{\sqrt{|g|}}, \quad
  M=\frac{h_{12}}{\sqrt{|g|}}, \quad
  N=\frac{h_{22}}{\sqrt{|g|}},
\quad |g|=det(g_{ij})=g_{11}g_{22}-g_{12}^2,
\end{equation}
$K(x,y)$ is the Gauss curvature that is determined by the relation:
\begin{equation}\label{5.22}
 K(x,y)=\frac{R_{1212}}{|g|},
\end{equation}
and
$$
R_{ijkl}=g_{lm}\big(\del_k\G^{(m)}_{ij}-\del_j\G^{(m)}_{ik}
  +\G^{(n)}_{ij}\G^{(m)}_{nk}-\G^{(n)}_{ik}\G^{(m)}_{nj}\big)
$$
is the Riemann curvature tensor depending on $g_{ij}$ and its
first and second derivatives.

Therefore, given a positive definite $(g_{ij})$, the Gauss-Codazzi
system consists of the three equations for the three unknowns $(L,
M, N)$ determining the second fundamental form. Note that, while
$(g_{ij})$ is positive definite, $R_{1212}$ may have any sign and
hence $K$ may have any sign.

From the viewpoint of geometry, the Gauss equation is a
Monge-Amp\`{e}re equation and the Codazzi equations are as
integrability relations. On the other hand, we are interested in a
fluid mechanic formulation for the isometric embedding problem so
that the problem may be solved via the approaches that have shown to
be useful in fluid mechanics, especially for the nonlinear partial
differential equations of mixed hyperbolic-elliptic type. To achieve
this, the way to reformulate the problem via solvability of the
Codazzi equations  under the Gauss equation has been
adopted in Chen-Slemrod-Wang \cite{CSW2}.

Set
$$
L=\rho v^2+p,  \quad M=-\rho uv, \quad  N=\rho u^2+p,
$$
choose the ``pressure'' $p$ as for the Chaplygin type gas: $
p=-\frac{1}{\rho}$,  and set $q^2=u^2+v^2$. Then the Codazzi
equations become the balance equations of momentum:
\begin{equation}\label{5.23}
\begin{split}
&\del_x(\rho uv)+\del_y(\rho v^2+p)
 =-\G^{(2)}_{22}(\rho v^2+p)-2\G^{(2)}_{12}\rho uv-\G^{(2)}_{11}(\rho u^2+p), \\
&\del_x(\rho u^2+p)+\del_y(\rho uv)
 =-\G^{(1)}_{22}(\rho v^2+p)-2\G^{(1)}_{12}\rho uv-\G^{(1)}_{11}(\rho u^2+p),
\end{split}
\end{equation}
and the Gauss equation becomes the Bernoulli law:
\begin{equation}\label{5.24-a}
\rho=\frac{1}{\sqrt{q^2+K}}.
\end{equation}

Set the ``sound" speed: $c=\sqrt{p'(\rho)}$, i.e.,
$c^2=\frac{1}{\rho^2}=q^2+K$. Then

\medskip
\begin{itemize}

\item When $K>0$, $\, q^2<c^2$ and the ``flow" is subsonic;

\item When $K<0$, $\, q^2>c^2$ and the ``flow" is supersonic;

\item When $K=0$, $\ q^2=c^2$ and the ``flow" is sonic.
\end{itemize}

\begin{figure}[h]
\centering
\includegraphics[scale=0.4]{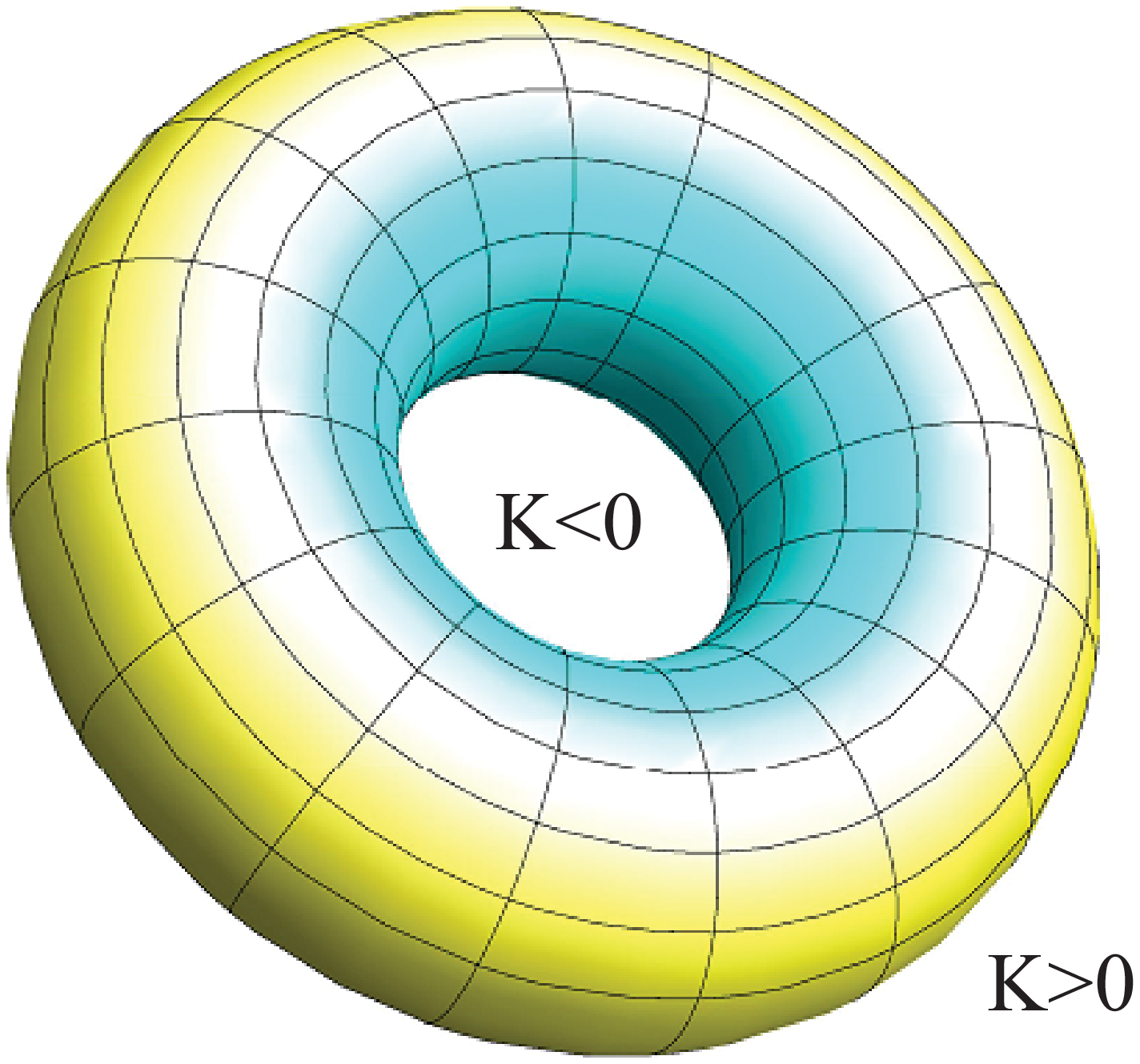}
\caption[]{The Gaussian curvature on a torus: Doughnut surface or
toroidal shell} \label{fig:Donut}
\end{figure}

Therefore, the existence of an isometric immersion is equivalent to
the existence of a weak solution of the balance equations of
momentum in the fluid mechanic formulation, which are nonlinear
partial differential equations of mixed hyperbolic-elliptic type.
Many usual surfaces have their Gauss curvature of changing sign,
such as tori: Doughnut surfaces or toroidal shells in Fig.
\ref{fig:Donut}.

An appropriate approximate method to construct approximate solutions
to the Gauss-Codazzi system has been designed, and  a compensated
compactness approach to establish the existence of a weak solution
has been developed in \cite{CSW2}. The advantage of the compensated
compactness approach is that it works for both the elliptic and
hyperbolic phase. For more details about isometric embedding of
Riemann manifolds in Euclidean spaces, we refer the reader to
Han-Hong \cite{HanHong} and  Yau \cite{Yau1}.

\smallskip
For isometric embedding of a higher dimensional
Riemannian manifold into the Euclidean space $\R^N$ with optimal
dimension $N$ (i.e., the Cartan-Janet dimension), the Gauss-Codazzi
equations should be supplemented by the Ricci equations to describe
the connection form on the normal bundle. The Gauss-Codazzi-Ricci
system even has no type in general, although it inherits important
geometric features, including the beautiful Div-Curl structure which
yields its weak continuity (see \cite{CSW3}). Also see
Bryant-Griffiths-Yang \cite{BGY}.

\subsection{Free-Boundary Techniques}

To explain the techniques, we focus on
the shock reflection-diffraction problem for potential flow which is
widely used in aerodynamics.

When a plane shock in the $(t,\x)$-coordinates,
$\x=(x_1,x_2)\in\R^2$, with the left state $(\rho,
\nabla_\x\Phi)=(\rho_1, u_1, 0)$ and the right state $(\rho_0, 0,0),
u_1>0, \rho_0<\rho_1$, hits a symmetric wedge
$$
W:=\{\x\, :\, |x_2|< x_1 \tan\theta_w, x_1>0\}
$$
head on, it experiences a reflection-diffraction process, where
$\rho$ is the density and $\Phi$ is the velocity potential of the
fluid; see Fig. \ref{fig:SRDP}. Then a self-similar reflected shock moves outward
as the original shock moves forward in time. The complexity of
reflection-diffraction configurations was first reported by Ernst
Mach in 1878, and experimental, computational, and asymptotic
analysis has shown that various different patterns may occur,
including regular and Mach reflection (cf.
\cite{BD,GlimmMajda,LaxLiu,Morawetz2,Serre,Neumann}). However, most
of the fundamental issues for shock reflection-diffraction have not
been understood, including the global structure, stability, and
transition of different patterns of shock reflection-diffraction
configurations. Therefore, it becomes essential to establish the
global existence and structural stability of solutions in order to
understand fully the shock reflection-diffraction phenomena.
Furthermore, this problem is also fundamental in the mathematical
theory of multidimensional conservation laws since their solutions
are building blocks and asymptotic attractors of general solutions
to the two-dimensional Euler equations for compressible fluid flow
(cf. Courant-Friedrichs \cite{CFr}, von Neumann \cite{Neumann}, and
Glimm-Majda \cite{GlimmMajda}; also see
\cite{BD,GQChen,Da,LaxLiu,Morawetz2,Serre}). As we will show below,
the problem involves nonlinear partial differential equations of
mixed elliptic-hyperbolic type, along with the other challenging
difficulties such as free boundary problems and corner
singularities, especially when a free boundary meets an elliptic
degenerate curve.

\begin{figure}[h]
\centering
\includegraphics[scale=0.9]{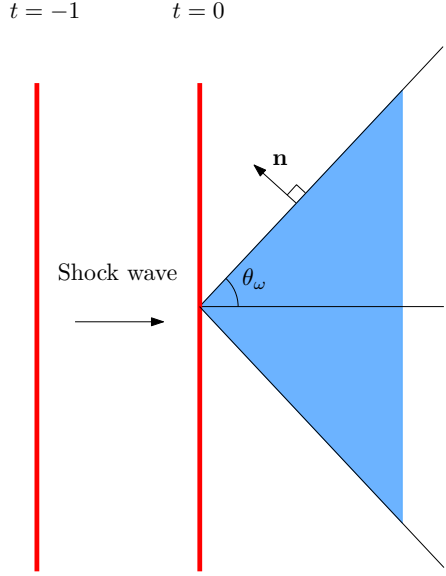}
\caption[]{Shock reflection-diffraction problem} \label{fig:SRDP-a}
\end{figure}

By scaling, the Euler equations for potential flow consist of the
conservation law of mass and the Bernoulli law with the form:
\begin{equation}\label{1.1.1}
\begin{cases}
\partial_t\rho + \nabla\cdot(\rho\nabla\Phi)=0,\\
\partial_t\Phi
+\frac{1}{2}|\nabla\Phi|^2+\frac{\rho^{\gamma-1}}{\gamma-1}
=\frac{\rho_0^{\gamma-1}}{\gamma-1},
\end{cases}
\end{equation}
where $\gamma>1$ is the adiabatic exponent of the
fluid. Then the reflection-diffraction problem can be formulated as
the following mathematical problem.

\medskip
{\bf Initial-Boundary Value Problem}. {\it Seek a solution of the
system with the initial condition at $t=0$:
\begin{equation}\label{initial-condition}
(\rho,\Phi)|_{t=0} =\begin{cases}
(\rho_0, 0) \qquad&  \mbox{for}\,\, |x_2|>x_1\tan\theta_w,\,\, x_1>0,\\
(\rho_1, u_1 x_1) \qquad &\mbox{for}\,\, x_1<0,
\end{cases}
\end{equation}
and the slip boundary condition along the wedge boundary $\partial
W$:
\begin{equation}\label{boundary-condition}
\nabla\Phi\cdot {\bf n}|_{\partial W}=0,
\end{equation}
where ${\bf n}$ is the exterior unit normal to $\partial W$. }

\medskip
Since the initial-boundary value problem is invariant under the
self-similar scaling:
$$
(t,{\x})\to (\alpha t, \alpha {\x}), \quad (\rho, \Phi)\to (\rho,
\frac{\Phi}{\alpha}) \qquad\quad \mbox{for $\alpha\ne 0$,}
$$
we seek self-similar solutions with the form:
$$
\rho(t, \x)=\rho(\xi,\eta), \quad \Phi(t, \x)=t\,\psi(\xi,\eta)
\qquad\quad \mbox{for}\quad (\xi,\eta)=\frac{{\bf x}}{t}.
$$
Then the pseudo-potential function
$\varphi=\psi-\frac{1}{2}(\xi^2+\eta^2)$ is governed by the
following potential flow equation of second-order:
\begin{equation}\label{5.24}
\dv\, \big(\rho(|D\varphi|^2, \varphi)D\varphi\big)
+2\rho(|D\varphi|^2, \varphi)=0
\end{equation}
with $\rho(|D\varphi|^2, \varphi)
:=\big(\rho_0^{\gamma-1}-(\gamma+1)(\varphi+\frac{1}{2}|D\varphi|^2)\big)^{\frac{1}{\gamma-1}},
$ where $D=(\p_\xi, \p_\eta)$. The sonic speed is
\begin{equation}\label{5.25}
c=c(|D\varphi|^2,\varphi,\rho_0^{\gamma-1})
:=\big(\rho_0^{\gamma-1}-(\gamma-1)(\varphi+\frac{1}{2}|D\varphi|^2)\big)^{1/2}.
\end{equation}
Equation \eqref{5.24} is a {\it second-order nonlinear partial
differential equations of mixed elliptic-hyperbolic type}. It is
strictly {\it elliptic} (i.e., pseudo-subsonic) if
\begin{equation}\label{5.26}
|D\varphi| < c(|D\varphi|^2,\varphi,\rho_0^{\gamma-1}),
\end{equation}
and strictly {\it hyperbolic} (i.e., pseudo-supersonic) if $|D\varphi|>
c(|D\varphi|^2,\varphi,\rho_0^{\gamma-1})$.

Shocks are discontinuities in the pseudo-velocity $D\varphi$. That
is, if $\Omega^+$ and $\Omega^-:=\Omega\setminus\overline{\Omega^+}$
are two nonempty open subsets of $\Omega\subset\R^2$ and
$S:=\partial\Omega^+\cap\Omega$ is a $C^1$ curve where $D\varphi$
has a jump,  then $\varphi\in W^{1,1}_{loc}(\Omega)\cap
C^1(\Omega^\pm\cup S)\cap C^2(\Omega^\pm)$ is a global weak solution
of \eqref{5.24} in $\Omega$ in the sense of distributions if and
only if $\varphi$ is in $W^{1,\infty}_{loc}(\Omega)$ and satisfies
equation \eqref{5.24} in $\Omega^\pm$ and the Rankine-Hugoniot
condition on $S$:
\begin{equation}\label{5.27}
\rho(|D\varphi|^2,\varphi)D\varphi\cdot {\bf n}\Big|_{S+}
=\rho(|D\varphi|^2,\varphi)D\varphi\cdot {\bf n}\Big|_{S-}.
\end{equation}
The continuity of $\varphi$ is followed by the continuity of the
tangential derivative of $\varphi$ across $S$. The discontinuity $S$
of $D\varphi$ is called a shock if $\varphi$ further satisfies the
physical entropy condition that the corresponding density function
$\rho(|D\varphi|^2,\varphi)$ increases across $S$ in the pseudo-flow
direction.

\medskip
The plane incident shock solution in the $(t, \x)$--coordinates
corresponds to a continuous weak solution $\varphi$ in the
self-similar coordinates $(\xi,\eta)$ with the following form:
\begin{eqnarray}
&&\varphi_0(\xi,\eta)=-\frac{1}{2}(\xi^2+\eta^2) \qquad
 \hbox{for } \,\, \xi>\xi_0,
 \label{flatOrthSelfSimShock1} \\
&&\varphi_1(\xi,\eta)=-\frac{1}{2}(\xi^2+\eta^2)+ u_1(\xi-\xi_0)
\qquad
 \hbox{for } \,\, \xi<\xi_0,
 \label{flatOrthSelfSimShock2}
\end{eqnarray}
respectively, where $\xi_0>0$ is the location of the incident shock,
uniquely determined by $(\rho_0,\rho_1,\gamma)$. Since the problem
is symmetric with respect to the axis $\eta=0$, it suffices to
consider the problem in the half-plane $\eta>0$ outside the
half-wedge $ \Lambda:=\{\xi\le 0,\eta>0\}\cup\{\eta>\xi
\tan\theta_w,\, \xi>0\}$. Then the initial-boundary value problem
can be formulated as the boundary value problem in $\Lambda$ in the
coordinates $(\xi,\eta)$.

\medskip
{\bf Boundary Value Problem}. {\it Seek a solution $\varphi$ of
equation \eqref{5.24} in the unbounded domain $\Lambda$ with the
slip boundary condition on the wedge boundary $\partial\Lambda$:
\begin{equation}\label{boundary-condition-3}
D\varphi\cdot {\bf n}|_{\partial\Lambda}=0
\end{equation}
and the asymptotic boundary condition at infinity:
\begin{equation}\label{boundary-condition-2}
\varphi\to \hat{\varphi}:=
\begin{cases} \varphi_0 \quad\mbox{for}\,\,\,
                         \xi>\xi_0, \eta>\xi \tan\theta_w,\\
              \varphi_1 \quad \mbox{for}\,\,\,
                          \xi<\xi_0, \eta>0,
\end{cases}
\mbox{when $\xi^2+\eta^2\to \infty$},
\end{equation}
in the sense that $ \displaystyle
\lim_{R\to\infty}\|\varphi-\hat{\varphi}\|_{C(\Lambda\setminus
B_R(0))}=0. $}

\begin{figure}[h]
\centering
\includegraphics[scale=0.8]{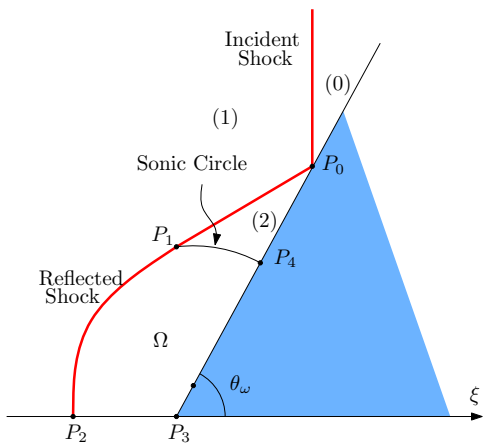}
\caption[]{Regular reflection-diffraction configuration}
\label{fig:RegularReflection}
\end{figure}

Since $\varphi_1$ does not satisfy the slip boundary condition, the
solution must differ from $\varphi_1$ in $\{\xi<\xi_0\}\cap\Lambda$,
and thus a shock diffraction by the wedge occurs. A local existence
theory of regular shock reflection-diffraction near the reflection
point $P_0$ can be established by following the von Neumann
detachment criterion \cite{Neumann} with the structure of solution
as in Fig. \ref{fig:RegularReflection}, when the wedge angle is
large and close to $\frac{\pi}{2}$, in which the vertical line is
the incident shock $S=\{\xi=\xi_0\}$ that hits the wedge at the
point $P_0=(\xi_0, \xi_0 \tan\theta_w)$, and state (0) and state (1)
ahead of and behind $S$ are given by $\varphi_0$ and $\varphi_1$,
respectively. The solutions $\varphi$ and $\varphi_1$ differ only in
the domain $P_0\PtUpL \PtLwL \PtLwR$ because of shock reflection by
the wedge boundary at $P_0$ and diffraction by the wedge vertex
$P_3$, where the curve $P_0\PtUpL \PtLwL$ is the reflected shock
with the straight segment $P_0\PtUpL$. State (2) behind $P_0\PtUpL$
can be computed explicitly with the form:
\begin{equation}\label{5.28}
\varphi_2(\xi,\eta)=-\frac{1}{2}(\xi^2+\eta^2)+u_2(\xi-\xi_0)+
(\eta-\xi_0\tan\theta_w)u_2\tan\theta_w,
\end{equation}
which satisfies $D\varphi\cdot {\bf n}=0$ on $\partial\Lambda\cap
\{\xi>0\}$;  the constant velocity $u_2$ and the slope of
$P_0\PtUpL$ are determined by $(\theta_w,\rho_0,\rho_1,\gamma)$ from
the two algebraic equations expressing \eqref{5.27} and continuous
matching of states (1) and (2) across $P_0\PtUpL$, whose existence
is exactly guaranteed by the condition on
$(\theta_w,\rho_0,\rho_1,\gamma)$ which is necessary for the existence
of a global regular shock reflection-diffraction configuration.

\medskip
A rigorous mathematical approach has been developed in Chen-Feldman
\cite{ChenFeldman06} to extend the local theory to a global theory
for solutions of regular shock reflection-diffraction, which
converge to the unique solution of the normal shock reflection when
$\theta_w$ tends to $\frac{\pi}{2}$. The solution $\varphi$ is
pseudo-subsonic within the sonic circle for state (2) with center
$(u_2, u_2\tan\theta_w)$ and radius $c_2>0$ (the sonic speed) and is
pseudo-supersonic outside this circle containing the arc
$\PtUpL\PtUpR$ in Fig. \ref{fig:RegularReflection}, so that
$\varphi_2$ is the unique solution in the domain $P_0\PtUpL \PtUpR$.
In the domain $\Omega$, the solution is pseudo-subsonic, smooth, and
$C^{1,1}$-matching with state (2) across $\PtUpL\PtUpR$ and
satisfies $\varphi_\eta=0$ on $\PtLwL\PtLwR$; the transonic shock
curve $\PtUpL\PtLwL$ matches up to second-order with $P_0\PtUpL$
continuously and is orthogonal to the $\xi$-axis at the point
$\PtLwL$ so that the standard reflection about the $\xi$--axis
yields a global solution in the whole plane.

One of the main difficulties for the global existence is that the
ellipticity condition for \eqref{5.24} is hard to control. The
second difficulty is that the ellipticity degenerates at the sonic
circle $\PtUpL\PtUpR$. The third difficulty is that, on
$\PtUpL\PtUpR$, the solution in $\Omega$ has to be matched with
$\varphi_2$ at least in $C^1$, i.e., the two conditions on the fixed
boundary $\PtUpL\PtUpR$: the Dirichlet and conormal conditions,
which are generically overdetermined for an elliptic equation since
the conditions on the other parts of boundary are prescribed. Thus
it is required to prove that, if $\varphi$ satisfies \eqref{5.24} in
$\Omega$, the Dirichlet continuity condition on the sonic circle,
and the appropriate conditions on the other parts of
$\partial\Omega$ derived from the boundary value problem, then the
normal derivative $D\varphi\cdot {\bf n}$ automatically matches with
$D\varphi_2\cdot {\bf n}$ along $\PtUpL\PtUpR$. Indeed, equation
\eqref{5.24}, written in terms of the function $w=\varphi-\varphi_2$ in the
$(x,y)$--coordinates defined near $\PtUpL\PtUpR$ such that
$\PtUpL\PtUpR$ becomes a segment on $\{x=0\}$, has the form:
\begin{equation}\label{5.29}
\big(2x-(\gamma+1)\partial_x w\big)\partial_{xx} w
+\frac{1}{c_2^2}\partial_{yy} w-\partial_x w=0
 \qquad\,\, \mbox{in } x>0 \mbox{  and near } x=0,
\end{equation}
plus the ``small'' terms that are controlled by
$\frac{\pi}{2}-\theta_w$ in appropriate norms. Equation \eqref{5.29}
is {\it strictly elliptic} if $\partial_x w<\frac{2x}{\gamma+1}$ and become
{\it degenerate elliptic} when $\partial_x w=\frac{2x}{\gamma+1}$. Hence, it
is required to obtain the $C^{1,1}$--estimates near $\PtUpL\PtUpR$
to ensure $|\partial_x w|<\frac{2x}{\gamma+1}$ which in turn implies both the
ellipticity of the equation in $\Omega$ and the match of normal
derivatives $D\varphi\cdot {\bf n}=D\varphi_2\cdot {\bf n}$ along
$\PtUpL\PtUpR$. Taking into account the ``small'' terms to be added
to equation \eqref{5.29}, it is needed to make the stronger estimate
$|\partial_x w|\le \frac{4x}{3(\gamma+1)}$ and assume that
$\frac{\pi}{2}-\theta_w$ is suitably small to control these
additional terms. Another issue is the non-variational structure and
nonlinearity of this problem which makes it hard to apply directly
the approaches of Caffarelli \cite{Ca} and Alt-Caffarelli-Friedman
\cite{AC,ACF}. Moreover, the elliptic degeneracy and geometry of the
problem makes it difficult to apply the hodograph transform approach
in Kinderlehrer-Nirenberg \cite{KinderlehrerNirenberg} and
Chen-Feldman \cite{ChenFeldman4} to fix the free boundary.

For these reasons, one of the new ingredients in the approach is to
develop further the iteration scheme in \cite{ChenFeldman1} to a
partially modified equation. Equation  \eqref{5.29} is modified in
$\Omega$ by a proper Shiffmanization (i.e. cutoff) that depends on
the distance to the sonic circle, so that the original and modified
equations coincide for $\varphi$ satisfying $|\partial_x w| \le
\frac{4x}{3(\gamma+1)}$, and the modified equation ${\mathcal
N}\varphi=0$ is elliptic in $\Omega$ with elliptic degeneracy on
$\PtUpL\PtUpR$. Then a free boundary problem is solved for this
modified equation: The free boundary is the curve $\PtUpL\PtLwL$,
and the free boundary conditions on $\PtUpL\PtLwL$ are
$\varphi=\varphi_1$ and the Rankine-Hugoniot condition \eqref{5.27}.
Moreover, the precise gradient estimate:
$|\partial_x w|<\frac{4x}{3(\gamma+1)}$ is made to ensure that $\varphi$
satisfies the original equation \eqref{5.24}.

\smallskip
This global theory for large-angle wedges has been extended in
Chen-Feldman \cite{ChenFeldman07} up to the sonic angle $\theta_s\le
\theta_c$, i.e. state (2) is sonic when $\theta_w=\theta_s$, such
that, as long as $\theta_w\in (\theta_s, \frac{\pi}{2}]$, the global
regular shock reflection-diffraction configuration exists, which
solves the von Neumann's sonic conjecture (1943) \cite{Neumann}
when $u_1<c_1$.
Furthermore, the optimal regularity of regular
reflection-diffraction solutions near the pseudo-sonic circle has
been shown in Bae-Chen-Feldman \cite{BCF} to be $C^{1,1}$ as established in
\cite{ChenFeldman06,ChenFeldman07}.

\smallskip
Some important efforts were made mathematically for the global
reflection problem via simplified models, including the unsteady
TSD equation and other nonlinear
models (cf. \cite{Keyfitz}). Furthermore, in order to deal with the
reflection problem, some asymptotic methods have been also developed
including the work by Lighthill, Keller, Morawetz, and others. The
physicality of weak Prandtl-Meyer reflection for supersonic
potential flow around a ramp has been also analyzed in
\cite{EllingLiu}. Another recent effort has been made on various
important physical, mixed elliptic-hyperbolic problems in steady
potential flow for which great progress has been made (cf.
\cite{GQChen,Chen4,ChenFeldman07,
morawetz-04} and
the references cited
therein).

\smallskip
The self-similar solutions for the full Euler equations are required
for the Mach reflection-diffraction configurations and general
two-dimensional Riemann problems. For this case, the self-similar
solutions are governed by a system that consists of two transport
equations and two {\it nonlinear equations of mixed
hyperbolic-elliptic type}. One of the important features in the
reflection-diffraction configurations is that the Euler equations
for potential flow and the full Euler equations coincide in some
important regions of the solutions.

\section{Singular Limits to Nonlinear Degenerate Hyperbolic Equations}

One of the important singular limit problems is the vanishing
viscosity limit problem for the Navier-Stokes equations to the Euler
equations for compressible barotropic fluids. The
Navier-Stokes equations for a compressible viscous, barotropic fluid
in Eulerian coordinates in $\mathbb{R}_+^2$ take the following form:
\begin{equation}
\label{Eq:NS-1}
\begin{cases}
\del_t\rho+ \del_x(\rho v)={}0,\\
\del_t(\rho v)+\del_x(\rho v^2{}+{}p)={}\e \del_{xx}v{},\\
\end{cases}
\end{equation}
with the initial conditions:
\begin{equation}
\label{Eq:Initial-Conditions} \rho(0,x){}={}\rho_0(x),\qquad
v(0,x){}={}v_0(x)
\end{equation}
such that $\lim_{x\to\pm\infty}(\rho_0(x), v_0(x))=(\rho^\pm,
v^\pm)$, where $\rho$ denotes the density, $v$ represents the
velocity of the fluid when $\rho>0$, $p$ is the pressure, $m=\rho v$
is the momentum, and $(\rho^\pm, v^\pm)$ are constant states with
$\rho^\pm>0$. The physical viscosity coefficient $\e$ is restricted
to $\e\in (0, \e_0]$ for some fixed $\e_0>0$.
For a polytropic perfect gases, the pressure-density relation is
determined by \eqref{3.15a}--\eqref{Eq:Pressure} with adiabatic
exponent $\gamma>1$.

Formally, when $\ve\to 0$, the Navier-Stokes equations become the
isentropic Euler equations \eqref{4.1a}, which is strictly
hyperbolic when $\rho>0$. However, near the vacuum $\rho=0$, the two
characteristic speeds of \eqref{4.1a} may coincide and the system be
degenerate hyperbolic.

The vanishing artificial/numerical viscosity limit to the isentropic
Euler equations with general $L^\infty$ initial data has been
studied by DiPerna \cite{DiPerna1}, Chen \cite{Chen1,Chen2a}, Ding
\cite{Ding}, Ding-Chen-Luo \cite{DCL}, Lions-Perthame-Souganidis
\cite{LPS}, and Lions-Perthame-Tadmor \cite{LPT}  via the methods of
compensated compactness. Also see DiPerna \cite{DiPerna2}, Morawetz
\cite{Mor}, Perthame-Tzavaras \cite{PerthameTzavaras}, and Serre
\cite{Serre-g} for the vanishing artificial/numerical viscosity
limit to general $2\times 2$ strictly hyperbolic systems of
conservation laws. The vanishing artificial viscosity limit to
general strictly hyperbolic systems of conservation laws with
general small $BV$ initial data was first established by
Bianchini-Bressan \cite{BB} via direct BV estimates with small
oscillation. Also see LeFloch-Westdickenberg \cite{LW} for the
existence of finite-energy solutions to the isentropic Euler
equations with finite-energy initial data for the case $1<\gamma\le
5/3$.

The idea of regarding inviscid gases as viscous gases with vanishing
real physical viscosity can date back the seminal paper by Stokes
\cite{Stokes} and the important contribution of Rankine
\cite{Rankine}, Hugoniot \cite{Hugoniot}, Rayleigh
\cite{Rayleigh}, and Taylor \cite{Taylor} (cf. Dafermos \cite{Da}).
However, the first
rigorous convergence analysis of vanishing physical viscosity from
the Navier-Stokes equations \eqref{Eq:NS-1} to the isentropic Euler
equations was made by Gilbarg \cite{Gi} in 1951, when he established
the mathematical existence and vanishing viscous limit of the
Navier-Stokes shock layers. For the convergence analysis confined in
the framework of piecewise smooth solutions; see Hoff-Liu \cite{HL},
G\`{u}es-M\'{e}tivier-Williams-Zumbrun \cite{GMWZ}, and the
references cited therein. The convergence of vanishing physical
viscosity with general initial data was first studied by
Serre-Shearer \cite{SS} for a $2\times 2$ system in nonlinear
elasticity with severe growth conditions on the nonlinear function
in the system.

In Chen-Perepelitsa \cite{CPer}, we have first developed new uniform
estimates with respect to the real physical viscosity coefficient
for the solutions of the Navier-Stokes equations with the
finite-energy initial data and established the $H^{-1}$-compactness
of weak entropy dissipation measures of the solutions of the
Navier-Stokes equations for any weak entropy pairs generated by
compactly supported $C^2$ test functions. With these, the existence
of measure-valued solutions with possibly unbounded support has been
established, which are confined by the Tartar-Murat commutator
relation with respect to two pairs of weak entropy kernels. Then we
have established the reduction of measure-valued solutions with
unbounded support for the case $\gamma\ge 3$ and, as corollary, we
have obtained the existence of global finite-energy entropy
solutions of the Euler equations with general initial data for
$\gamma\ge 3$. We have further simplified the reduction proof of
measure-valued solutions with unbounded support for the case
$1<\gamma\le 5/3$ in LeFloch-Westdickenberg \cite{LW} and extended
to the whole interval $1<\gamma<3$ . With all of these, we have
established the first convergence result for the vanishing physical
viscosity limit of solutions of the Navier-Stokes equations to a
finite-energy entropy solution of the isentropic Euler equations
with general finite-energy initial data. We remark that
the
existence of finite-energy solutions to the isentropic Euler
equations
for the case $\gamma>5/3$ have been
established, which is in addition to the result in \cite{LW} for
$1<\gamma\le 5/3$.

More precisely, consider the Cauchy problem
\eqref{Eq:NS-1}--\eqref{Eq:Initial-Conditions} for the Navier-Stokes
equations in $\R_+^2$. Hoff's theorem in \cite{Hoff} (also see Kanel
\cite{Kanel} for the case of the same end-states) indicates that,
when the initial functions $(\rho_0(x), v_0(x))$ are smooth with the
lower bounded density $\rho_0(x)\ge c^\e_0>0$ for $x\in\R$ and
$$
\lim_{x\to\pm\infty}(\rho_0(x), v_0(x))=(\rho^\pm, v^\pm),
$$
then there exists a unique smooth solution $(\rho^\e(t,x),
v^\e(t,x))$, globally in time, with $\rho^\e(t,x)\ge c_\e(t)$ for
some $c_\e(t)>0$ for $t\ge 0$ and
$\lim_{x\to\pm\infty}(\rho^\e(t,x), v^\e(t,x))=(\rho^\pm, v^\pm)$.

\medskip
Let $(\bar{\rho}(x),\bar{v}(x))$ be a pair of smooth monotone
functions satisfying $(\bar{\rho}(x), \bar{v}(x))=(\rho^\pm, v^\pm)$
when $\pm x\ge L_0$ for some large $L_0>0$.
The total mechanical
energy for \eqref{Eq:NS-1} in $\R$ with respect to the pair
$(\bar{\rho},\bar{v})$ is
\begin{equation}
\label{def:energy} E[\rho, v](t):=\int_{\R}\Big(\eta^*(\rho,
m)-\eta^*(\bar{\rho}, \bar{m})-\nabla\eta^*(\bar{\rho},
\bar{m})\cdot (\rho-\bar{\rho}, m-\bar{m})\Big)\, dx\ge 0,
\end{equation}
where $\bar{m}=\bar{\rho}\bar{v}$.

Combining the uniform estimates and the compactness of weak entropy
dissipation measures in $H^{-1}_{loc}$ with the compensated
compactness argument  and the reduction of the measure-valued
solution $\nu_{t,x}$, we conclude

\begin{theorem}[Chen-Perepelitsa \cite{CPer}]\label{main-theorem}
Let the initial functions $(\rho_0^\e, v_0^\e)$ be smooth and
satisfy the following conditions: There exist $E_0, E_1, M_0>0$,
independent of $\e$, and $c_0^\e>0$ such that
\begin{enumerate}
\item[(i)] $\rho_0^\e(x)\ge c_0^\e>0, \quad \int\rho_0^\e(x)|v_0^\e(x)-\bar{v}(x)|\,dx \le M_0
<\infty$;

\medskip
\item[(ii)] The total mechanical energy with respect to
$(\bar{\rho}, \bar{u})$ is finite:
$$
\int\Big(\frac{1}{2}\rho_0^\e(x)|v_0^\e(x)-\bar{v}(x)|^2
 +h(\rho_0^\e(x),\bar{\rho}(x))\Big)dx\le E_0<\infty,
$$
where $h(\rho,
\bar{\rho})=\rho e(\rho)-\bar{\rho}e(\bar{\rho})-
\big(e(\bar{\rho})+\bar{\rho}e'(\bar{\rho})\big)(\rho-\bar{\rho})$;
\item[(iii)]
$\e^2\int \frac{|\rho_{0,x}^\e(x)|^2}{\rho_0^\e(x)^3} \,dx\le
E_1<\infty$;

\item[(iv)]
$(\rho^\e_0(x), \rho^\e_0(x)v^\e_0(x))\to (\rho_0(x),
\rho_0(x)v_0(x))$ in the sense of distributions as $\e\to 0$, with
$\rho_0(x)\ge 0$ a.e.
\end{enumerate}
Let $(\rho^\e, m^\e),
m^\e=\rho^\e v^\e$, be the solution of the Cauchy problem
\eqref{Eq:NS-1}--\eqref{Eq:Initial-Conditions} for the Navier-Stokes
equations with initial data $(\rho^\e_0(x), \rho^\e_0(x)v^\e_0(x))$ for each
fixed $\e>0$. Then, when $\e\to 0$, there exists a subsequence of
$(\rho^\e, m^\e)$ that converges almost everywhere to a
finite-energy entropy solution $(\rho, m)$ of the Cauchy problem
\eqref{4.1a} and \eqref{Eq:Initial-Conditions} with initial data
$(\rho_0(x), \rho_0(x)v_0(x))$ for the isentropic Euler equations
with $\gamma>1$.
\end{theorem}

\bigskip
{\bf Acknowledgments}. Gui-Qiang Chen's research was supported in part
by the National Science Foundation under Grants DMS-0935967,
DMS-0807551, and DMS-0505473, the Natural Science Foundation of
China under Grant NSFC-10728101, and the Royal Society-Wolfson
Research Merit Award (UK). The author would like to thank his
collaborators M. Bae, C.~M. Dafermos, M. Feldman, K.
Karlsen, Ph. LeFloch, M. Perepelitsa, B. Perthame, M. Slemrod, D.
Wang, among others, for their explicit/implicit contributions
reflected here. This paper was written as part of the
International Research Program on Nonlinear Partial Differential
Equations at the Centre for Advanced Study at the Norwegian Academy
of Science and Letters in Oslo during the Academic Year 2008--09.

\end{document}